\documentclass[11pt,reqno]{amsart}
\usepackage[T1]{fontenc}
\usepackage{lmodern}
\usepackage{amssymb,mathtools,microtype}
\usepackage[letterpaper,margin=1.12in,headsep=13pt]{geometry}
\usepackage[hidelinks]{hyperref}
\makeatletter
\renewcommand{\@defaultbiblabelstyle}[1]{[#1]}
\makeatother
\numberwithin{equation}{section}
\newtheorem{theorem}{Theorem}[section]
\newtheorem{lemma}[theorem]{Lemma}
\newtheorem{proposition}[theorem]{Proposition}
\theoremstyle{remark}
\newtheorem{remark}[theorem]{Remark}
\theoremstyle{plain}
\hypersetup{pdftitle={On the higher-order Titchmarsh divisor problem},
pdfauthor={Guangwei Hu, Guangshi Lu, Zhiwei Wang, Jie Wu},
pdfsubject={Higher-order Titchmarsh divisor correlations}}

\begin{document}
\title[On the higher-order Titchmarsh divisor problem]
{On the higher-order Titchmarsh divisor problem}
\author[G. Hu]{Guangwei Hu}
\address[Guangwei Hu]{School of Mathematics and Statistics, Shandong Normal University, Jinan 250100, Shandong, China}
\email{gwhu@sdnu.edu.cn}
\author[G. L\"u]{Guangshi L\"u}
\address[Guangshi L\"u]{School of Mathematics, Shandong University, Jinan 250100, Shandong, China}
\email{gslv@sdu.edu.cn}
\author[Z. Wang]{Zhiwei Wang}
\address[Zhiwei Wang]{School of Mathematics, Shandong University, Jinan 250100, Shandong, China}
\email{zhiwei.wang@sdu.edu.cn}
\author[J. Wu]{Jie Wu}
\address[Jie Wu]{Mathematical Research Center, Shandong University, Jinan 250100, Shandong, China}
\email{jie.wu@sdu.edu.cn}

\keywords{Titchmarsh divisor problem, shifted primes, divisor functions, sieve methods, dilation
graphs.}
\subjclass[2020]{11N37; 11N36; 11N13}
\date{\today}

\begin{abstract}
For every fixed integer $k\geq2$ and every fixed nonzero integer $h$, we prove that, as $x\to\infty$,
\begin{equation*}
\sum_{x<n\leq2x}\Lambda(n)\tau_k(n+h)
\sim C_{k,h}x(\log x)^{k-1},
\end{equation*}
where $\tau_k$ is the $k$-fold divisor function and $C_{k,h}>0$ is a constant depending only on $k$ and $h$.
\end{abstract}
\maketitle

\section{Introduction}

In number theory, a general hypothesis asserts that the multiplicative properties of the integers $n$ and $n+h$    are independent for any $h\in\mathbb{Z}\setminus\{0\}$.  In general, problems merging these two structures, such as the celebrated twin prime conjecture and $abc$-conjecture, are very difficult. In this paper, we are concerned with the Titchmarsh divisor problem, which is also in this class.

For an integer $k\geq2$, we denote by $\tau_k$ the $k$-fold divisor function. The higher-order Titchmarsh divisor problem concerns the average of $\tau_k(p+h)$ over primes $p$, where $h$ is a fixed nonzero integer. For any fixed integers $k\ge2$ and $h\ne0$, one may expect the asymptotic formula 
\begin{equation}\label{asymp}
	\sum_{x<n\le2x}\Lambda(n)\tau_k(n+h)
	\sim C_{k,h}x(\log x)^{k-1}
	\qquad (x\to\infty),
\end{equation}
where  $\Lambda$ denotes the von Mangoldt function and $C_{k,h}>0$ is a constant depending only on $k$ and $h$.
We refer to, e.g., Matom\"aki--Radziwi{\l}{\l}--Tao
\cite[Conjecture~1.1(iii) and Remark 1.2]{Matomaki-Radziwill-Tao-2019a}
for the asymptotic \eqref{asymp} and local density factor $C_{k,h}$ in the higher order Titchmarsh divisor problem.

For $k=2$, the function $\tau_2=\tau$ is the usual divisor function, and the problem reduces to the classical Titchmarsh
divisor problem.  Titchmarsh \cite{Titchmarsh-1930} established an asymptotic formula under the generalized Riemann hypothesis (GRH) for Dirichlet $L$-functions. For the shift $-h$ with $h>0$ fixed, his formula reads
\begin{equation}\label{eq:classical}
\sum_{h<p\leq x}\tau(p-h)=\frac{\varphi(h)}{h}\prod_{p\nmid h}\Big(1+\frac{1}{p(p-1)}\Big)x+O_h\Big(\frac{x\log\log x}{\log x}\Big),
\end{equation}
where $\varphi$ is Euler's totient function. In order to prove this, Titchmarsh opened the divisor function. Thus Dirichlet's hyperbola method transformed the problem to count primes in arithmetic progressions with moduli up to $\sqrt{x}$. The GRH provides the required estimates for primes in these arithmetic progressions. 30 years later, Linnik \cite{Linnik-1963} used his dispersion method \cite{Linnik-1960} to remove the GRH. Subsequently some simpler proofs have been found thanks to the Bombieri--Vinogradov theorem; see \cite{Halberstam-1967,Rodriquez-1965}.

In 1982, Deshouillers--Iwaniec \cite{Deshouillers-Iwaniec-1982} developed spectral estimates for sums of Kloosterman sums.
These estimates provide cancellation that can be used to study primes in arithmetic progressions to moduli slightly beyond $x^{1/2}$.
Combining these new methods with combinatorial identities,
Fouvry \cite{Fouvry-1985} and Bombieri--Friedlander--Iwaniec \cite{Bombieri-Friedlander-Iwaniec-1986} independently
improved the error term in \eqref{eq:classical} (with a slightly more precise main term) to $O_{h, A}(x (\log x)^{-A})$ for an arbitrary number $A>0$. It is also interesting to study the uniformity with respect to the shift $h$, see e.g. the work of  Fiorilli \cite{Fiorilli-2012}, Drappeau \cite{Drappeau-2017}, Drappeau--Topacogullari \cite{Drappeau-Topacogullari-2019}, Assing--Blomer--Li \cite{Assing-Blomer-Li-2021}.

For $k\geq3$, the asymptotic formula \eqref{asymp}
has not yet been established unconditionally. The standard sieve argument would
give an asymptotic formula with error $o(x(\log x)^{k-1})$ if
$\Lambda$ had a level of distribution exceeding $1-1/k$.
The Elliott--Halberstam conjecture would give this range, but
it remains open.

It is also of interest to study the problem on average over the shift.
Matom\"aki--Radziwi{\l}{\l}--Tao
\cite[Theorem~1.3(iii)]{Matomaki-Radziwill-Tao-2019a} proved that,
for each fixed integer $k\geq2$ and fixed constants $A>0$
and $0<\varepsilon<1/2$,
\begin{align*}
\sum_{x<n\leq2x}\Lambda(n)\tau_k(n+h)=xP_{k,h}(\log x)+O_{k,A,\varepsilon}\Big(x(\log x)^{-A}\Big)
\end{align*}
for all but $O_{k,A,\varepsilon}(H(\log x)^{-A})$ integers $h$ with $0<|h|\leq H$, whenever $x^{8/33+\varepsilon}\leq H\leq x^{1-\varepsilon}$. Here $P_{k,h}$ is a computable polynomial of degree $k-1$.
Their subsequent paper \cite[Theorem~1.1(ii)]{Matomaki-Radziwill-Tao-2019b} allows $H$ to be as small as $(\log X)^{10000k\log k}$,
with a weaker error term. Further developments include asymptotic formulae for multipoint correlations of primes and of divisor functions along arithmetic progressions, valid for almost all common differences in suitable ranges \cite{Matomaki-Radziwill-Shao-Tao-Teravainen-2026}. 

In this paper, we establish an  asymptotic formula for every fixed divisor order and every fixed nonzero shift.

\begin{theorem}\label{thm:main}
Let $k\geq2$ and $h\neq0$ be fixed integers. As $x\to\infty$, we have
\begin{equation*}
\sum_{x<n\leq2x}\Lambda(n)\tau_k(n+h)=C_{k,h}x(\log x)^{k-1}+o_{k,h}\Big(x(\log x)^{k-1}\Big),
\end{equation*}
where
\begin{equation}\label{eq:constant}
	C_{k,h}=\frac{1}{(k-1)!}\prod_{p\mid h}(1-p^{-1})^{k-1}\prod_{p\nmid h}\frac{1-p^{-1}(1-p^{-1})^{k-1}}{1-p^{-1}}.
\end{equation}
\end{theorem}

When $k=2$, the local factor at $p\nmid h$ is $1+1/(p(p-1))$, so the theorem recovers the classical Titchmarsh asymptotic formula. For each fixed $k$, all local factors are positive, and those at $p\nmid h$ equal $1+O_k(p^{-2})$.  For $k\geq3$, Theorem~\ref{thm:main} establishes an asymptotic formula for every fixed nonzero shift, with the predicted arithmetic constant.

\subsection*{Notation}

Throughout the proof, $k\geq2$ and $h\neq0$ are fixed integers, and $L=\log x$. The letter $p$ denotes a prime.
We write $\tau=\tau_2$ and use $\mu$ for the M\"obius function. We write $e(t)=\exp(2\pi it)$ and $\mathbb T=\mathbb R/\mathbb Z$.
For $t\in\mathbb R$, the notation $\|t\|_\mathbb T$ denotes the distance from $t$ to the nearest integer.
For a positive integer $n$, the symbols $v_p(n)$ and $P^-(n)$ denote its $p$-adic valuation and smallest prime factor, respectively, with $P^-(1)=\infty$. An integer is $z$-sifted if $P^-(n)>z$. Unless stated otherwise, $\Phi\in C_c^\infty((0,\infty))$ is fixed and nonnegative. Implied constants may depend on $k,h,\Phi$ and
on the auxiliary parameters fixed before $x\to\infty$.
Every arithmetic argument is required to be positive; in particular,
we restrict $v=n-h$ to $v\geq1$.

\section{The density model and auxiliary weights}

\subsection{The density model and moment bounds}
We shall introduce the density model, as in \cite[Section~2]{Walsh-2026}. For a finite set of primes $\mathcal U$, we denote by $\mathcal U'$ the set of all primes outside $\mathcal U$. For each positive integer $n$, we define
\begin{align*}
n_{\mathcal U}=\prod_{p\in\mathcal U}p^{v_p(n)}.
\end{align*}
Thus the factor $n_{\mathcal U}$ contains the full power of every prime in $\mathcal U$ that divides $n$, and $n=n_{\mathcal U}n_{\mathcal U'}$. A positive integer $a$ is $\mathcal U$-smooth if all its prime factors lie in $\mathcal U$, or equivalently if $a=a_{\mathcal U}$. We write $(n,\mathcal U)=1$ if no prime in $\mathcal U$ divides $n$. 
 
 We shall use the corresponding notations
\begin{align*}
y=\exp(L^{2/5}),\quad \mathcal S=\{p:p\leq y\},\quad V(t)=\prod_{p\leq t}(1-p^{-1}),\quad V=V(y),\quad E=\exp(\sqrt L).
\end{align*}
For each positive integer $n$, the factors $n_{\mathcal S}$ and $n_{\mathcal S'}$ are coprime by definition. The multiplicativity
of $\tau_k$ therefore gives
\begin{align*}
\tau_k(n)=\tau_k(n_{\mathcal S})\tau_k(n_{\mathcal S'}).
\end{align*}
Next, we replace the second by an approximation to the mean of $\tau_k$ over integers whose prime factors all exceed $y$. The corresponding Dirichlet series is
\begin{align*}
\sum_{\substack{n\ge1\\(n,\mathcal S)=1}}\frac{\tau_k(n)}{n^s}=\zeta(s)^k\prod_{p\le y}(1-p^{-s})^k\quad(\Re s>1).
\end{align*}
The pole at $s=1$ gives the model density. For $t\ge x^{1/2}$, the density is
\begin{equation}\label{eq:model-density}
M_k(t)=\operatorname*{Res}_{s=1}\Big\{\zeta(s)^k\prod_{p\le y}(1-p^{-s})^k t^{s-1}\Big\}.
\end{equation}
The residue is a polynomial of degree $k-1$ in $\log t$. By the bounds in \cite[(2.4)--(2.6)]{Walsh-2026}, we have
\begin{align*}
M_k(t)=\frac{V^k(\log t)^{k-1}}{(k-1)!}\Big(1+O_k(L^{-3/5})\Big)
\end{align*}
uniformly for $x^{1/2}\le t\le x^2$. The same bounds give positivity on $[x^{1/2}/2,x^2]$ for sufficiently large $x$, since $\log t\asymp L$. We obtain a smooth nonnegative extension to $0<t<x^{1/2}$ by multiplying this expression by a fixed nonnegative smooth cutoff in $t/x^{1/2}$, equal to one for $t\ge x^{1/2}$ and to zero for $t\le x^{1/2}/2$.  The coefficient bounds and the choice of cutoff imply that, for every fixed integer $r\ge0$,
\begin{align*}
\big|(t\partial_t)^rM_k(t)\big|\ll_{k,r}V^kL^{k-1}\quad(0<t\le x^2).
\end{align*}
With this extension, we define the density model for $\tau_k$ by
\begin{align*}
D_k(n)=\tau_k(n_{\mathcal S}) \frac{M_k(n_{\mathcal S'})}{V}.
\end{align*}
Since $V$ is the density of integers with no prime factor at most $y$, the quotient $M_k/V$ models the average divisor weight among these integers. The factor $\tau_k(n_{\mathcal S})$ preserves the exact contribution of the $y$-smooth part of $n$.
We denote the remainder by
\begin{align*}
\xi_k(n)=\tau_k(n)-D_k(n).
\end{align*}
For estimates in absolute value, we use the nonnegative auxiliary function
\begin{align*}
Q_k(n)=\tau_k(n)+D_k(n),
\end{align*}
which satisfies $|\xi_k(n)|\le Q_k(n)$.

For positive averages, we shall also use the simpler model
\begin{align*}
D_k^0(n)=B_k\tau_k(n_{\mathcal S}){ \text{ with }} B_k=\frac{V^{k-1}L^{k-1}}{(k-1)!}.
\end{align*}
The next lemma compares $D_k$ and $D_k^0$ for $n\asymp x$ with $n_{\mathcal S}\le E$.


\begin{lemma}\label{lem:model}
For $n\asymp x$ with $n_{\mathcal S}\leq E$, we have
\begin{equation}\label{eq:model}
D_k(n)=D_k^0(n)\Big(1+O(L^{-1/2})\Big).
\end{equation}
For every positive integer $n\leq x^2$, we also have $0\leq D_k(n)\ll L^{k-1}\tau_k(n)$ and $D_k(n)\ll D_k^0(n)$.
\end{lemma}

\begin{proof}
Since $n\asymp x$ and $n_{\mathcal S}\leq E$, it follows from \cite[Lemma~2.1]{Walsh-2026} that
\begin{align*}
M_k(n_{\mathcal S'})/V=B_k(1+O_k(L^{-1/2})).
\end{align*}
The definition of $D_k$ then yields \eqref{eq:model}. The uniform bound for $M_k$ gives
\begin{align*}
0\leq D_k(n)\ll_k V^{k-1}L^{k-1}\tau_k(n_{\mathcal S})=(k-1)!\,D_k^0(n).
\end{align*}
The remaining estimate follows from $\tau_k(n_{\mathcal S})\le\tau_k(n)$ and $V\le1$.
\end{proof}

The two moment estimates below are \cite[(2.13)--(2.14), Lemma~2.2]{Walsh-2026}. We also record their consequence for a fixed product of shifted divisor weights.

\begin{lemma}\label{lem:moments}
The integers $a,r\geq1$ and the constants $c>0$ and
$0<\delta\leq C$ are fixed. Uniformly for
$x^\delta\leq N\leq x^C$, we have
\begin{equation}\label{eq:moments}
\sum_{n\leq cN}\tau_a(n)^r n_{\mathcal S}^{rL^{-2/5}}
\ll NL^{a^r-1}.
\end{equation}
Consequently,
\begin{equation}\label{eq:tail}
\sum_{n\leq cN}\tau_a(n)^r\mathbf 1_{n_{\mathcal S}>E}
\ll NL^{a^r-1}\exp(-rL^{1/10}).
\end{equation}
\end{lemma}

\begin{proof}
Both estimates follow from \cite[Lemma~2.2]{Walsh-2026},
applied in the wider range $x^{\delta/2}\leq N\leq x^C$.
This also covers the case $\delta=C$.
\end{proof}

We may also multiply the summand in \eqref{eq:tail} by a fixed
finite product of shifted divisor functions. Over a fixed
interval $n\asymp x$, the resulting sum is bounded by
$xL^{O(1)}\exp(-c'L^{1/10})$, where $c'>0$ depends only
on the fixed parameters.

Indeed, we write the additional weight as
$\prod_{j=1}^{m}\tau_{b_j}(n+h_j)$.
By H\"older's inequality, we derive
\begin{align*}
\sum_{n\asymp x}\tau_a(n)^r
\mathbf 1_{n_{\mathcal S}>E}
\prod_{j=1}^{m}\tau_{b_j}(n+h_j)\leq
\Big(\sum_{n\asymp x}\tau_a(n)^{r(m+1)}
\mathbf 1_{n_{\mathcal S}>E}\Big)^{\frac{1}{m+1}}
\prod_{j=1}^{m}
\Big(\sum_{n\asymp x}\tau_{b_j}(n+h_j)^{m+1}\Big)^{\frac{1}{m+1}}.
\end{align*}
We estimate the first factor by \eqref{eq:tail}, with $r$
replaced by $r(m+1)$. For each remaining factor, we translate
the summation variable and apply \eqref{eq:moments} without
the factor involving $n_{\mathcal S}$.
The resulting bound is
$xL^{O(1)}\exp(-rL^{1/10})$, so we may take $c'=r$.

\subsection{Auxiliary weights}\label{sec:marks}

We introduce auxiliary weights depending on prime divisors
in fixed intervals, as in \cite[Section~2.3]{Lu-Wang-Wu-2026}.
We choose an integer $K\geq1$ and exponents
$1/10<a_1<\cdots<a_K<1/5$. For sufficiently large $x$,
we define the pairwise disjoint sets of primes
\begin{align*}
\mathcal B_i=\{p:\exp(L^{a_i})\leq p\leq\exp(2L^{a_i})\}.
\end{align*}
The integer $K$ and the exponents $a_i$ are fixed before
$x\to\infty$. We set $\mathcal B=\bigcup_{i=1}^K\mathcal B_i$ and
\begin{align*}
\omega_i(n)=\sum_{p\in\mathcal B_i}\mathbf 1_{p\mid n},
\quad V_i=\sum_{p\in\mathcal B_i}\frac1p=\log2+o(1).
\end{align*}
For $0<q<1$, we define
\begin{equation}\label{eq:marks}
\mathcal W_q(n)=\prod_{i=1}^K
\frac{q^{\omega_i(n)-1}\omega_i(n)}{V_i},
\qquad \rho=\frac1{2k},\quad \mathcal W=\mathcal W_\rho.
\end{equation}
The asymptotic for $V_i$ follows from Mertens' theorem.
We assign the value zero to a factor with $\omega_i(n)=0$.
The weight $\mathcal W(n)$ satisfies
\begin{align*}
0\leq\mathcal W(n)\leq
\Big(\frac{2}{\log2}\sup_{j\geq1}j\rho^{j-1}\Big)^K<\infty
\end{align*}
for sufficiently large $x$. The local mean of $\mathcal W$ involves the factors
\begin{equation}\label{eq:local-mass}
J_p=\int_{\mathbb Z_p} \tau_k(p^{v_p(t)})\mathbf 1_{p\nmid t-h}\,dt=\begin{cases}
(1-p^{-1})^{1-k}-p^{-1}, & p\nmid h,\\
1-p^{-1}, & p\mid h.
\end{cases}
\end{equation}
Here $v_p$ also denotes the valuation on $\mathbb Z_p$, and
$dt$ is Haar measure normalized by $\int_{\mathbb Z_p}dt=1$.
For each $\mathcal S$-smooth integer $a$, the normalized
local-density weight is
\begin{align*}
\mu_*(a)=\prod_{p\leq y}\frac1{J_p}
\int_{\substack{t\in\mathbb Z_p\\v_p(t)=v_p(a)}}
\tau_k(p^{v_p(a)})\mathbf 1_{p\nmid t-h}\,dt.
\end{align*}
The definition of $J_p$ gives $\sum_{a=a_{\mathcal S}}\mu_*(a)=1$.
For a bounded function $F$ that depends only on the $y$-smooth
part of its argument, we write
\begin{align*}
\langle F\rangle_*=
\sum_{a=a_{\mathcal S}}\mu_*(a)F(a),
\qquad m_{\mathcal W}=\langle\mathcal W\rangle_*.
\end{align*}
Sums of functions depending on disjoint sets of prime valuations
factor as products of the corresponding local sums.
For sufficiently large $x$, every prime in the sets $\mathcal B_i$
is coprime to $h$.  {The symbol $\vartheta_p$ denotes the sum of
$\mu_*(a)$ over $\mathcal S$-smooth $a$ divisible by $p$.
The local densities give}
\begin{align*}
\vartheta_p
=\frac{(1-p^{-1})^{1-k}-(1-p^{-1})}{J_p}
=\frac{k}{p}+O_k(p^{-2}).
\end{align*}
For each set, we distinguish one prime divisor of the argument
in the factor $\omega_i(n)\rho^{\omega_i(n)-1}$.  {The factorization of the local-density weights gives}
\begin{align*}
m_{\mathcal W}
=\prod_{i=1}^K\Big\{
\frac{1}{V_i}\sum_{p\in\mathcal B_i}\vartheta_p
\prod_{\substack{q\in\mathcal B_i\\q\ne p}}
\Big(1-(1-\rho)\vartheta_q\Big)\Big\}.
\end{align*}
Since $\sum_{p\in\mathcal B_i}p^{-2}=o(1)$ and
$V_i=\log2+o(1)$, the $i$th factor is
$k\exp(k(\rho-1)V_i)+o(1)$. Hence, for fixed $K$,
\begin{align*}
m_{\mathcal W}=\Big(k2^{k(\rho-1)}\Big)^K+o(1).
\end{align*}
Since $0\leq\mathcal W(n)\ll_{K,k}1$, we have
$\mathcal W(n)^2\ll_{K,k}\mathcal W(n)$. Hence there are
positive constants $c_{K,k}$ and $C_{K,k}$ such that
\begin{equation}\label{eq:mark-moments}
0<c_{K,k}\leq m_{\mathcal W}\leq C_{K,k},\quad
 {\langle(\mathcal W/m_{\mathcal W})^2\rangle_*}\leq C_{K,k}.
\end{equation}
These constants are independent of the fixed exponents $a_i$;
the threshold for $x$ may depend on them.

\section{Arithmetic progressions and the normalized sieve}

We begin with an estimate for $\xi_k$ in arithmetic progressions.
We then impose the condition that $n-h$ be sifted, first up to $y$
and then up to a small fixed power of $x$. Throughout this section,
$\mathcal W$ has the form \eqref{eq:marks}, with the number of
prime sets and their exponents fixed before $x\to\infty$.
All sums involving $n-h$ are restricted to $n-h\geq1$.

\subsection{An estimate over arithmetic progressions}

Throughout the proof, $c_0=1/(2(k+2))$. We first establish a progression estimate that allows an
auxiliary weight. Since $\mathcal W(n)$ depends only on
$n_{\mathcal S}$, we state the result for a general weight $H$
with this property.

\begin{lemma}\label{lem:progressions}
For every function $H$ satisfying $H(n)=H(n_{\mathcal S})$
and $|H(n)|\leq L^C$, where $C$ is fixed, and every fixed $A>0$,
we have
\begin{align*}
  \sum_{r\le x^{c_0}}\Big| \sum_{r\mid n-h}\Phi(n/x)H(n)\xi_k(n)\Big|\ll_{A,k,h,C,\Phi}xL^{-A}.
\end{align*}
The estimate is uniform over all such functions $H$, including
functions that depend on $x$.
\end{lemma}

\begin{proof}
We first remove the terms with $n_{\mathcal S}>E$.
After we interchange the sums over $r$ and $n$, their contribution is
\begin{align*}
\ll_\Phi L^C\sum_{n\asymp x}\Big(\tau_k(n)+D_k(n)\Big)\tau(n-h)\mathbf 1_{n_{\mathcal S}>E}.
\end{align*}
Since $D_k(n)\ll_k L^{k-1}\tau_k(n)$, we apply Lemma~\ref{lem:moments} with the additional weight $\tau(n-h)$ to bound this contribution by
\begin{align*}
O\Big(xL^{O_{k,C}(1)}\exp(-c_kL^{1/10})\Big).
\end{align*}
We fix $a=n_{\mathcal S}\leq E$ and write $n=at$. Then
$(t,\mathcal S)=1$ and
\begin{equation}
  \xi_k(at)=\tau_k(a)\Big(\tau_k(t)-M_k(t)/V\Big).
  \label{eq:centering}
\end{equation}
With $g=(a,r)$, the congruence $at\equiv h\; (\bmod r)$ has no solution unless $g\mid h$. Otherwise, with $q=r/g$, it is equivalent to
\begin{align*}
  t\equiv b\; (\bmod q),\quad b\equiv (h/g)(a/g)^{-1}\; (\bmod q).
\end{align*}
If $(b,q)>1$, a prime dividing $h$ must divide $t$. Such a prime belongs to $\mathcal S$ for large $x$, so there is again no contribution. For each fixed pair $(a,q)$, there are at most $\tau(|h|)$ possible values of $r$. Therefore, we may work with reduced residue classes and put $Y=x/a=x^{1-o(1)}$.

For a fixed integer $j\geq1$, a scale $ {Y_0}>0$ and a reduced
residue class $b\;(\bmod q)$, it follows from Mellin inversion and character orthogonality that
\begin{equation}\label{eq:divisor-ap}
\begin{aligned}
\sum_{m\equiv b\;(\bmod q)}\tau_j(m)\psi(m/ {Y_0})=\frac{1}{\varphi(q)}\mathop{\mathrm{Res}}_{s=1}
\Big\{\zeta(s)^j\prod_{p\mid q}(1-p^{-s})^j
 {Y_0}^s\widehat\psi(s)\Big\}+O( {Y_0}^{1/2}q^{j/2}),
\end{aligned}
\end{equation}
where $\psi$ is a smooth function supported in a fixed compact subinterval of $(0,\infty)$, and
$\widehat\psi(s)=\int_0^\infty\psi(t)t^{s-1}\,dt$.  Indeed, we shift the contour in each character integral to $\Re s=1/2$. The convexity bound gives a conductor factor $q^{j/4+\varepsilon}$. We choose $\varepsilon<j/4$ and use the rapid decay of $\widehat\psi$ to obtain the factor $q^{j/2}$; see \cite[Chapter~5]{Iwaniec-Kowalski-2004}. The argument requires only finitely many derivatives of $\psi$. If these derivatives are bounded by fixed powers of $L$, we incur at most an additional factor $L^{O(1)}$ in the error term.

We write $t=t_1\cdots t_k$ as an ordered factorization.
The condition that $t$ be $y$-sifted excludes
 {the divisibility conditions} $p\mid t_i$, where $p\le y$, $p\nmid q$ and $1\le i\le k$.  {We index these conditions by the pairs} $(p,i)$ and apply the Bonferroni inequalities up to degree $R=\lceil C_A\log L\rceil$.  {For each collection of these conditions, we define}
\begin{align*}
a_i=\prod_{(p,i)\,\mathrm{selected}}p \quad(1\le i\le k).
\end{align*}
 {These conditions impose} $a_i\mid t_i$ for each $i$, and
\begin{align*}
  a_1\cdots a_k\le y^{R+1} =\exp(O_A(L^{2/5}\log L))=x^{o(1)}.
\end{align*}
We write $t_i=a_iv_i$ and apply \eqref{eq:divisor-ap} to the
remaining sum.  {The number of such collections, including}
those in the first omitted degree, is $x^{o(1)}$.
The absolute Bonferroni error is bounded by the nonnegative sum of degree $R+1$. Its reciprocal coefficient mass is at most
\begin{align*}
\frac{\Big(k\sum_{p\leq y}p^{-1}\Big)^{R+1}}{(R+1)!}.
\end{align*}
We apply \eqref{eq:divisor-ap} to this nonnegative sum, using a fixed nonnegative smooth majorant of the cutoff if necessary. The residue term is bounded by
\begin{align*}
  \frac{Y}{\varphi(q)}L^{O_k(1)} \frac{\Big(k\sum_{p\leq y}p^{-1}\Big)^{R+1}}{(R+1)!}.
\end{align*}
At degree $j$, the derivatives needed to recover the full Euler
product from the truncated residue introduce an additional factor
$j^{k-1}$. Since $\sum_{p\leq y}p^{-1}\ll\log L$, we may
choose $C_A$ sufficiently large and apply Stirling's formula
to obtain any prescribed negative power of $L$.
The main term of the divisor sum over $y$-sifted integers is
\begin{equation}\label{eq:rough-residue}
  \frac{1}{\varphi(q)}\int_0^\infty \Phi(t/Y)M_{k,q}(t)\,dt,
\end{equation}
where
\begin{align*}
G_{j,q}(s)=\prod_{\substack{p\mid q\\p>y}}(1-p^{-s})^j,\; M_{k,q}(t)=\mathop{\mathrm{Res}}_{s=1}
\Big\{\zeta(s)^k\prod_{p\leq y}(1-p^{-s})^k
G_{k,q}(s)t^{s-1}\Big\}.
\end{align*}
We treat the term $M_k(t)/V$ in \eqref{eq:centering} by applying the same argument with one colour and the cutoff
$\Phi(t/Y)M_k(t)/V$.  {The derivative bounds for $M_k$ in Section~2 show that
this cutoff has derivatives of every fixed order in $t/Y$
bounded by fixed powers of $L$.} We obtain the main term
\begin{equation}\label{eq:model-residue}
\frac{G_{1,q}(1)}{\varphi(q)}\int_0^\infty \Phi(t/Y)M_k(t)\,dt.
\end{equation}
Uniformly for $q\leq x^{c_0}$ and $0\leq i\leq k-1$, we have
\begin{align*}
  |G_{k,q}^{(i)}(1)-\mathbf 1_{i=0}|+|G_{1,q}(1)-1|\ll_k e^{-L^{2/5}}L^{O_k(1)}.
\end{align*}
We obtain this bound by differentiating the logarithm: $q$ has
at most $O(L)$ distinct prime divisors, and every prime in these
products exceeds $y$. Thus the difference between \eqref{eq:rough-residue} and \eqref{eq:model-residue} is
exponentially small, up to a fixed power of $L$ and the factor
$Y/\varphi(q)$.

We now sum the errors with coefficients $\tau_k(a)H(a)$. The errors
in the main terms are $O(xL^{-A})$ for any prescribed $A>0$, since
\begin{align*}
  \sum_{a=a_{\mathcal S}}\frac{\tau_k(a)}{a}=V^{-k}=L^{O_k(1)},
  \quad
  \sum_{q\le x^{c_0}}\frac{1}{\varphi(q)}\ll L.
\end{align*}
The truncated coefficients have total absolute mass $x^{o(1)}$.
Thus the total contour error from \eqref{eq:divisor-ap} is
\begin{align*}
  x^{1/2+c_0(1+k/2)+o(1)}=x^{3/4+o(1)}.
\end{align*}
This proves the lemma.
\end{proof}

\subsection{The condition \texorpdfstring{$P^-(n-h)>y$}{P-minus(n-h) > y}}

We next impose the condition $P^-(n-h)>y$. The Bonferroni
inequalities reduce this condition to the progression sums of
Lemma~\ref{lem:progressions}, together with a positive error term.
We also allow a divisor sum supported on $y$-sifted integers;
this will be needed for the next sieve step.

For $v=n-h\geq1$, we write $S_y(v)=\mathbf 1_{P^-(v)>y}$, and let $\omega(d)$ denote the number of distinct prime divisors of $d$.
For each integer $R\geq0$, we define
\begin{equation}\label{eq:bonferroni}
\begin{aligned}
T_R(v)=\sum_{\substack{d\mid\prod_{p\leq y}p\\\omega(d)\leq R}}
\mu(d)\mathbf 1_{d\mid v},\;
E_{R+1}(v)=\sum_{\substack{d\mid\prod_{p\leq y}p\\\omega(d)=R+1}}
\mathbf 1_{d\mid v}.
\end{aligned}
\end{equation}
By the Bonferroni inequalities, $|S_y(v)-T_R(v)|\leq E_{R+1}(v)$.

\begin{lemma}\label{lem:small-primes}
The coefficients $\lambda_r$ satisfy $|\lambda_r|\leq1$
and vanish unless $r\leq x^{c_0/2}$ and $P^-(r)>y$. We define
$A_\lambda(v)=\sum_{r\mid v}\lambda_r$ and
$A_{|\lambda|}(v)=\sum_{r\mid v}|\lambda_r|$.
For every fixed $A>0$, we have
\begin{equation}\label{eq:small-prime-sieve}
  \sum_n\Phi(n/x)\xi_k(n)\mathcal W(n)S_y(v)A_\lambda(v)\ll xL^{-A}.
\end{equation}
The estimate holds uniformly for all coefficients satisfying these conditions.
\end{lemma}

\begin{proof}
We choose $R=\lceil C_A\log L\rceil$.
Each $y$-smooth divisor $d$ in \eqref{eq:bonferroni} satisfies
$d\leq y^{R+1}=x^{o(1)}$.
The prime factors of $d$ are at most $y$, whereas those of $r$
are larger. The product $dr$ determines $d$ and $r$ uniquely. Thus Lemma~\ref{lem:progressions} gives
\begin{equation}
  \sum_n\Phi(n/x)\xi_k(n)\mathcal W(n)T_R(v)A_\lambda(v)\ll xL^{-A}.
  \label{eq:truncated-sieve}
\end{equation}
We apply the same argument to the first omitted degree and obtain
\begin{equation}
  \begin{aligned}
&\sum_n\Phi(n/x)\tau_k(n)\mathcal W(n)E_{R+1}(v)A_{|\lambda|}(v)\\
&=\sum_n\Phi(n/x)D_k(n)\mathcal W(n)
E_{R+1}(v)A_{|\lambda|}(v)+O(xL^{-A}).
  \end{aligned}
  \label{eq:positive-error}
\end{equation}
Both sums are nonnegative. Their sum bounds the error from
the replacement of $T_R$ by $S_y$, since $|\xi_k|\leq\tau_k+D_k$.
Since $E_{R+1}(v)$ and $A_{|\lambda|}(v)$ are each at most
$\tau(v)$, Lemma~\ref{lem:moments} bounds the contribution
from $n_{\mathcal S}>E$ to the sum with $D_k$ by
$xL^{O(1)}\exp(-cL^{1/10})$ for some $c>0$.
For $n_{\mathcal S}\leq E$, Lemma~\ref{lem:model} gives
$D_k\ll_k D_k^0$.

We fix $a=n_{\mathcal S}\le E$ and discard the restriction that $n/a$ has no prime factor at most $y$. The conditions $a\mid n$ and
$dr\mid n-h$ are compatible only if $(a,d)\mid h$, since $(a,r)=1$. In the compatible case, their density is $(a,d)/(adr)\le |h|/(adr)$,
and the smooth count has error $O_\Phi(1)$. The total counting error is $x^{c_0/2+o(1)}$. Consequently the positive sum with $D_k^0$
is at most
\begin{equation}
  \begin{aligned}
&\sum_n\Phi(n/x)D_k^0(n)\mathcal W(n)E_{R+1}(v)A_{|\lambda|}(v)\ll_{k,h,K,\Phi}
 xB_kV^{-k}L\,\frac{(\sum_{p\leq y}p^{-1})^{R+1}}{(R+1)!}
 +x^{c_0/2+o(1)}.
  \end{aligned}
  \label{eq:sieve-tail}
\end{equation}
Since $B_kV^{-k}L=L^k/((k-1)!V)$ is bounded by a fixed power of $L$, we choose $C_A$ sufficiently large and apply Stirling's
formula to bound the right-hand side of \eqref{eq:sieve-tail} by $O_{A,K}(xL^{-A})$. Together with \eqref{eq:positive-error},
we can replace $T_R$ by $S_y$ in \eqref{eq:truncated-sieve} with an error of the same order.
\end{proof}

\subsection{Progressions for the positive model}
 We evaluate the corresponding positive model
sum and retain the mean of $\mathcal W$ in its main term.

\begin{lemma}\label{lem:model-progressions}
For $P^-(r)>y$ and $r\le x^{c_0/2}$, we write
\begin{align*}
  \sum_{r\mid n-h}\Phi(n/x)D_k^0(n)\mathcal W(n)S_y(n-h)=\frac{x}{r}\Big(\int_0^\infty\Phi(t)\,dt\Big) B_km_{\mathcal W}\prod_{p\le y}J_p+ {\mathcal R}(r).
\end{align*}
For every fixed $A>0$, we have
\begin{align*}
  \sum_{\substack{r\le x^{c_0/2}\\P^-(r)>y}}| {\mathcal R}(r)|\ll xL^{-A}.
\end{align*}
\end{lemma}

\begin{proof}
We first restrict to integers $n$ whose $y$-smooth part $a=n_{\mathcal S}$ satisfies $a\leq E$. After summing over $r$, we bound the omitted contribution
by a fixed divisor moment with the additional weight $\tau(n-h)$. Lemma~\ref{lem:moments} then gives the bound $xL^{O(1)}\exp(-cL^{1/10})$ for some $c>0$. We may also restrict to $(a,h)=1$, since otherwise a prime $p\leq y$ divides both $n$ and $h$, and hence $S_y(n-h)=0$.

For $(a,h)=1$, we write $n=at$. The conditions $a=n_{\mathcal S}$ and $S_y(at-h)=1$ require $p\nmid t$ and $p\nmid at-h$ for every prime $p\leq y$. They exclude the residue class $0\;(\bmod p)$ and, when $p\nmid ah$, the distinct class $ha^{-1}\;(\bmod p)$.
The number of excluded classes is denoted by $\nu_a(p)$. Thus $\nu_a(p)=1$ when $p\mid ah$ and $\nu_a(p)=2$ otherwise.
Their combined density is
\begin{align*}
\delta_a=\prod_{p\leq y}(1-\nu_a(p)/p).
\end{align*}
If $\delta_a=0$, the local conditions have no solution, and the corresponding sum is zero. We omit these values of $a$
before normalization. In particular, when both $h$ and $a$ are odd, the conditions exclude both residue classes modulo $2$.

Since $(a,r)=1$, the congruence $at\equiv h\;(\bmod r)$
determines a unique residue class modulo $r$. Moreover, $P^-(r)>y$
ensures that $r$ is coprime to all the small-prime moduli.
The Chinese remainder theorem therefore combines these conditions
without any further compatibility restriction.
We apply the Bonferroni inequalities to the forbidden residue classes with $R=\lceil C_A\log L\rceil$. By the Chinese remainder
theorem, each compatible set of conditions determines a single residue class with modulus at most $ry^{R+1}$. We count each class with the smooth weight and incur an error $O_\Phi(1)$.
The retained terms, including the values $a\leq E$, have total
absolute coefficient mass $x^{o(1)}$. Thus, after summing over
$r\leq x^{c_0/2}$, we obtain a total counting error of
$x^{c_0/2+o(1)}$.

For the main terms, we bound the truncation error by
\begin{align*}
\frac{x}{ar}L^{O(1)}
\frac{(2\sum_{p\leq y}p^{-1})^{R+1}}{(R+1)!}.
\end{align*}
We multiply this bound by $B_k\tau_k(a)\mathcal W(a)$ and sum over
$a$ and $r$. These sums introduce at most a fixed power of $L$,
with a constant depending on the fixed choice of auxiliary prime sets.
For sufficiently large $C_A$, Stirling's formula then gives
a total error of $O_{A,K}(xL^{-A})$.

For each fixed $a$, the remaining main term is
\begin{align*}
  \frac{x}{ar}\Big(\int_0^\infty\Phi(t)\,dt\Big)
     B_k\tau_k(a)\mathcal W(a)\delta_a.
\end{align*}
We also bound the tail $a>E$ in the model sum.
With $\theta=1/\log y$, we apply Rankin's trick to obtain
the upper bound
\begin{align*}
  E^{-\theta}\|\mathcal W\|_\infty
    \prod_{p\le y}(1-p^{-1+\theta})^{-k}
  \ll_K L^{O_k(1)}e^{-L^{1/10}}.
\end{align*}
We may therefore remove the restriction $a\leq E$ from the model sum.
We evaluate the complete sum through its local factors. For $p\nmid h$,
the term with $v_p(a)=0$ has mass $1-2/p$, and the term with
$v_p(a)=j\ge 1$ has mass $(1-p^{-1})\tau_k(p^j)/p^j$. For $p\mid h$,
only $v_p(a)=0$ is allowed, with mass $1-p^{-1}$. The sums of these local masses are the constants $J_p$ in
\eqref{eq:local-mass}. Hence
\begin{align*}
  \sum_{\substack{a=a_{\mathcal S}\\(a,h)=1}}
      \frac{\tau_k(a)\mathcal W(a)\delta_a}{a}
    =m_{\mathcal W}\prod_{p\le y}J_p.
\end{align*}
We have retained the factor $1/r$ throughout the calculation.
The error estimates above remain valid after summation over $r$,
which proves the lemma.
\end{proof}

\subsection{The condition \texorpdfstring{$P^-(n-h)>x^\gamma$}{P-minus(n-h) > x-to-gamma}}
We now sieve the primes between $y$ and $x^\gamma$.
The following estimate retains the local mean $m_{\mathcal W}$
in its error term, so we can choose the number of auxiliary prime
sets after the required sieve accuracy.

\begin{proposition}\label{prop:presieve}
For a sufficiently large even integer $s$, we set
$\gamma=c_0/(8s+8)$ and $z=x^\gamma$. Then
\begin{equation}
  \limsup_{x\to\infty}
  \frac{\Big|\sum_n\Phi(n/x)\xi_k(n)\mathcal W(n)
                     \mathbf 1_{P^-(n-h)>z}\Big|}
       {xL^{k-2}m_{\mathcal W}}
    \ll_{k,h,\Phi}\frac{e^{-s}}{\gamma}.
  \label{eq:presieve}
\end{equation}
The implied constant is independent of $s$, $K$ and the fixed
choice of auxiliary prime sets. We take the limit with all these
parameters fixed.
\end{proposition}

\begin{proof}
We apply the block sieve of \cite[Lemma~2.9]{OpenAI-2026f}, in the form
given in \cite[Lemma~2.6]{OpenAI-2026c}, to the primes $y<p\leq z$,
with density $g(p)=1/p$. We obtain upper and lower sieve
polynomials
$\lambda^\pm(v)=\sum_{r\mid v}\lambda_r^\pm$,
whose coefficients satisfy $|\lambda_r^\pm|\leq1$ and vanish
unless $r$ is squarefree and
\begin{align*}
P^-(r)>y,\quad r\leq z^{4s+2}<x^{c_0/2}.
\end{align*}
These polynomials satisfy
\begin{equation}\label{eq:sieve-gap}
\begin{gathered}
\lambda^-(v)
\leq \mathbf 1_{(v,\prod_{y<p\leq z}p)=1}
\leq \lambda^+(v),\;
G(v):=\lambda^+(v)-\lambda^-(v)\geq0,\;
\sum_r\frac{\lambda_r^+-\lambda_r^-}{r}
\ll e^{-s}\frac{V(z)}{V(y)}.
\end{gathered}
\end{equation}
The implied constant is absolute, since $g(p)\leq1/2$ and
$\sum_{w<p\leq w^2}g(p)\ll1$ uniformly for $w>1$.
We apply the upper and lower sieve bounds separately to the
nonnegative weights $\Phi\mathcal W \tau_k S_y$ and
$\Phi\mathcal W D_k S_y$. By Lemma~\ref{lem:small-primes}, we obtain
\begin{equation}\label{eq:positive-sieve}
  \begin{aligned}
&\Big|\sum_n\Phi(n/x)\xi_k(n)\mathcal W(n)
 \mathbf 1_{P^-(n-h)>z}\Big|\\
&\leq\sum_n\Phi(n/x)D_k(n)\mathcal W(n)S_y(n-h)G(n-h)
 +O(xL^{-A}).
  \end{aligned}
\end{equation}
Since $G\geq0$, we apply Lemma~\ref{lem:model} to replace
$D_k$ by $D_k^0$ in the sum over $n_{\mathcal S}\leq E$,
with a relative error of $O_k(L^{-1/2})$.
For $n_{\mathcal S}>E$, we use $G(v)\leq2\tau(v)$
and Lemma~\ref{lem:moments} to bound the contribution by
$xL^{O(1)}\exp(-cL^{1/10})$ for some $c>0$. By Lemma~\ref{lem:model-progressions}, the remaining sum with $D_k^0$
equals
\begin{equation}\label{eq:model-sieve-gap}
  x\Big(\int_0^\infty\Phi(t)\,dt\Big)B_km_{\mathcal W}
    \prod_{p\leq y}J_p
    \sum_r\frac{\lambda_r^+-\lambda_r^-}{r}
  +O_{A,K}(xL^{-A}).
\end{equation}
We substitute this expression into \eqref{eq:positive-sieve},
keeping the factor $1+O_k(L^{-1/2})$ from the model comparison.
Once all auxiliary parameters are fixed, this relative error
is absorbed into the final little-$o$ term.

We extract $V(y)$ from $B_k\prod_{p\leq y}J_p$.
The remaining local factor is
$(1-p^{-1})^{k-2}J_p$. By \eqref{eq:local-mass}, it equals
$(1-p^{-1})^{k-1}$ if $p\mid h$, and
$(1-p^{-1})^{-1}-p^{-1}(1-p^{-1})^{k-2}$ otherwise. The latter factors are $1+O_k(p^{-2})$, so
\begin{equation}\label{eq:normalization}
  B_k\prod_{p\leq y}J_p=(C_{k,h}+o(1))L^{k-1}V(y).
\end{equation}
We use \eqref{eq:normalization}, \eqref{eq:sieve-gap} and
Mertens' estimate $V(z)\asymp1/(\gamma L)$.
They bound \eqref{eq:model-sieve-gap}, including the relative
model error, by
\begin{align*}
  C_{k,h,\Phi}\frac{e^{-s}}{\gamma}xL^{k-2}m_{\mathcal W}
  +o_K(xL^{k-2}).
\end{align*}
After division by $xL^{k-2}m_{\mathcal W}$, we let $x\to\infty$.
The remainders tend to zero because $m_{\mathcal W}\geq c_K>0$.
The implied constant is independent of
$\|\mathcal W\|_\infty/m_{\mathcal W}$.
\end{proof}

\section{Two correlation estimates}

We set $Z=\exp(L^{0.24})$.
The Type II estimate replaces the prime indicators in a convolution
by weights supported on $Z$-sifted integers.
We then prove a correlation estimate for the resulting convolution,
which allows auxiliary prime divisors at both endpoints.

\subsection{The Type II estimate}
The next lemma extends \cite[Theorem~3.1]{OpenAI-2026c} to
any fixed nonzero shift. The weights $\mathcal W_q$ are those
in \eqref{eq:marks}.

\begin{lemma}\label{lem:type-ii}
The parameters $\delta,C,A>0$ and $q\in(0,1)$ are fixed,
as is the nonzero integer $h$. The scales $M,N$ satisfy
$M,N\geq x^\delta$ and $MN\asymp x$. The function $F$ satisfies
$|F(n)|\leq1$ and $F(pn)=F(n)$ for every
positive integer $n$ and every prime $p\in\bigcup_i\mathcal B_i$, where
the sets $\mathcal B_i$ define $\mathcal W_q$. For an interval $I\subseteq[M,2M]$, we define the coefficients by
\begin{align*}
 \alpha_m=m^{it}\Big\{\mathbf 1_{m\ \mathrm{prime}}
 -\frac{\mathbf 1_{P^-(m)>Z}}{V(Z)\log m}\Big\},
 \quad m\in I,\quad |t|\le L^C,
\end{align*}
with $\alpha_m=0$ outside $I$. The coefficients $\beta_n$
satisfy $|\beta_n|\leq L^C$ and vanish unless
$N\leq n\leq2N$ and $P^-(n)>Z$.
If $K$ is sufficiently large in terms of the fixed parameters, then
\begin{equation}\label{eq:type-ii}
 \sum_{m,n}\alpha_m\beta_nF(mn+h)\mathcal W_q(mn+h)
 \ll MN L^{-A}.
\end{equation}
The lower bound for $K$ does not depend on the exponents
defining the sets $\mathcal B_i$.
\end{lemma}

\begin{proof}
The case $h=-1$ is \cite[Theorem~3.1]{OpenAI-2026c}.
For a general fixed $h$, we use $a=-h$ and $X=MN$.
The proof follows the expansion in \cite[(3.8)]{OpenAI-2026c}.
We check the determinant equation and the estimates on both
the minor and major arcs.

We first expand the product of the $\omega_i$ in $\mathcal W_q$
and select one prime from each $\mathcal B_i$. If a selected prime
$p$ divides $mn-a$ to exponent at least two, then
$mn\equiv a\;(\bmod p^2)$. For sufficiently large $x$, every
selected prime is coprime to $a$. If $p\mid m$, this congruence
has no solution; otherwise, it determines a unique residue
class for $n$ modulo $p^2$. Since $p^2=o(N)$ uniformly over
the selected primes, each such class contains $O(N/p^2)$
admissible values of $n$. We therefore bound the total
contribution by
\begin{align*}
XL^{O(1)}\sum_{p\in\bigcup_i\mathcal B_i}p^{-2},
\end{align*}
where the implied constant may depend on the fixed prime sets.

As in \cite[Section~3.1]{OpenAI-2026c}, we restrict the selected products
$A$ and $B$ to $Y\leq A,B\leq4Y$ and write
$mn-a=An_0$ and $rs-a=Bn_0$. Outside the prime-square terms already
discarded, we have $F(An_0)=F(n_0)$ and
$\omega_i(An_0)=1+\omega_i(n_0)$. Thus the factor depending on $n_0$
is $F(n_0)q^{\sum_i\omega_i(n_0)}$, whose absolute value is at most
one. We may restore the discarded terms after this factorization
at the same negligible cost. We denote the factored sum
with products of this size by $S_Y$. The Cauchy--Schwarz inequality gives
$|S_Y|^2\ll(X/Y)Q_Y$, where
\begin{align*}
 Q_Y={}&\Big(\prod_iV_i^{-2}\Big)
 \sum_{\substack{A,B,m,n,r,s\\
 mn-a=An_0,\ rs-a=Bn_0\ \text{for some }n_0}}
 \eta(A/Y)\eta(B/Y)\alpha_m\beta_n
 \overline{\alpha_r\beta_s}.
\end{align*}
Here $A$ and $B$ are products of one prime from each $\mathcal B_i$;
we retain the original dyadic supports of $m,n,r,s$.
The function $\eta$ is the smooth cutoff in the partition of unity of
\cite[Section~3.1]{OpenAI-2026c}.

For products sharing a selected prime, we use
\begin{equation}\label{eq:affine-fourth-divisor}
 \sum_{\substack{n_0\ll X/Y\\1\le a+\ell n_0\ll X}}
 \tau(a+\ell n_0)^2\ll_a(X/Y)L^3,
 \quad \ell\asymp Y,\quad(\ell,a)=1,
\end{equation}
where $\ell$ is a selected product. Indeed, $\tau(n)^2\le \tau_4(n)$. In each ordered factorization
of $n$ into four factors, we omit a largest factor and obtain
\begin{align*}
 \tau_4(n)\le4\sum_{\substack{d\mid n\\d\le n^{3/4}}}\tau_3(d).
\end{align*}
If $n=a+\ell n_0$, every divisor $d\mid n$ is coprime to $\ell$.
Thus $n_0$ lies in one residue class modulo $d$. We count the
integers in each class and sum over $d$ to obtain
\begin{align*}
 \sum_{\substack{n_0\ll X/Y\\1\le a+\ell n_0\ll X}}
 \tau(a+\ell n_0)^2
 \ll_a(X/Y)L^3+X^{3/4}L^2\ll_a(X/Y)L^3,
\end{align*}
since $Y=x^{o(1)}$. Together with the count in
\cite[Section~3.1]{OpenAI-2026c}, \eqref{eq:affine-fourth-divisor}
allows us to omit products with a common selected prime.
Their contribution to $Q_Y$ is $O(XY L^{-E_0})$ for every
fixed $E_0>0$.

When $(A,B)=1$, the two equations involving $n_0$ are equivalent to
\begin{align*}
T=Bmn-Ars=a(B-A).
\end{align*}
To see the converse, we reduce this equation modulo $A$ and $B$
and obtain $A\mid mn-a$ and $B\mid rs-a$. The two quotients
are equal and positive on the ranges under consideration.

In the construction of \cite[(3.17)]{OpenAI-2026c}, we use the vectors
$P=(DBm,s)$ and $Q=(DAr,n)$, where $D$ is the product of the
additional primes introduced there. Every prime factor of $AD$ or $BD$
is smaller than $Z$, whereas $m,n,r,s$ are $Z$-sifted. The preceding
discard also ensures $(A,B)=1$. A prime dividing both coordinates
of either vector must exceed $Z$ and divide $T$. Since
$|T|\ll_aY<Z$ on the kernel support, this would force $T=0$.
But a prime dividing $B$ divides none of $A,r,s$, so $Bmn-Ars\ne0$.
Thus both vectors are primitive.
We use the major arcs $\mathfrak M$ defined in
\cite[(3.11)]{OpenAI-2026c}. For the minor-arc argument, we choose a real even function
$\psi_a\in C_c^\infty(\mathbb R)$ equal to one on
$[-4|a|,4|a|]$. The factor in \cite[(3.15)]{OpenAI-2026c} becomes
\begin{equation}\label{eq:shift-kernel}
 \psi_a(T/Y)\Big(\mathbf 1_{T=a(B-A)}
 -\int_{\mathfrak M}e\Big(\theta(T-aB+aA)\Big)\,d\theta\Big).
\end{equation}
We retain the first admissibility condition on root ratios
from \cite[Section~3.2]{OpenAI-2026c}. In the second condition, we replace
the phases $DUT_fr_P$ by $aDUT_fr_P$. For a primitive integer
vector $P$ with $P_1\neq0$, we choose $w_P\in\mathbb Z^2$
such that $\det(P,w_P)=1$ and define
$r_P=(w_P)_1/P_1$ modulo one. Here $U$, denoted by $Z$ in the cited construction, is the
product of the other entries transferred between the prime lists, and $T_f$ denotes the product of
the newly selected primes.

For fixed newly selected primes, distinct products $DU$ and
$D'U'$ give the phase difference $a(DU-D'U')T_fr_P$.
Since $a(DU-D'U')T_f$ is a nonzero integer, multiplication
by this coefficient preserves Haar measure on ${\mathbb T}$. The same exponential saving therefore follows from the pair count
and {Chebyshev's inequality for the nonnegative counting function in}
\cite[Lemma~3.2]{OpenAI-2026c}.

We average over the root ratios by \cite[Lemma~3.4]{OpenAI-2026c}.
The substitution $r_P\mapsto ar_P$ increases the number of
interval endpoints in each admissibility condition by at most
a factor depending on $a$. The bound $|T|\ll_a Y$ changes
the bounds on the path coordinates by only a fixed factor.
The powers of $L$ in the path counts, congruence moduli and
boundary counts therefore remain unchanged.
For an edge with no prime stored or transferred from the additional
lists, we use the lattice blocks of \cite[Section~4.5]{OpenAI-2026c}. The new cutoff affects only the bounded number of neighbouring
cells and the constants in \cite[(4.27)]{OpenAI-2026c}. In the integer coordinates ${\mathbf z},{\mathbf z}'$ of
that formula, the kernel is
\begin{align*}
\begin{aligned}
 K_a(({\mathbf z},B),({\mathbf z}',A))=\frac1Y\psi_a\Big(\frac{\det({\mathbf z},{\mathbf z}')}{Y}\Big)\int_{\mathbb T\setminus\mathfrak M}
 e\Big(\theta(\det({\mathbf z},{\mathbf z}')-aB+aA)\Big)\,d\theta,
 \end{aligned}
\end{align*}
The kernel $K_{\mathrm{full}}$ has unrestricted integer
coordinates $b,c$ in place of $aB,aA$.
On the corresponding counting-measure spaces, we define
\begin{align*}
 (J_af)({\mathbf z},b)=
 \begin{cases}f({\mathbf z},b/a),&a\mid b,\\0,&a\nmid b.\end{cases}
\end{align*}
The map $J_a$ is an isometry for every nonzero integer $a$.
Consequently,
\begin{equation}\label{eq:kernel-compression}
 K_a=J_a^*K_{\mathrm{full}}J_a,
 \quad\|K_a\|\le\|K_{\mathrm{full}}\|.
\end{equation}
The endpoint multipliers in \cite[(4.28)]{OpenAI-2026c} contribute only
a factor depending on $K$. We take the Fourier transform
in the unrestricted coordinate and obtain the matrices in
\cite[(4.30)]{OpenAI-2026c}, indexed by $\theta\notin\mathfrak M$.
We now apply \cite[(4.31)--(4.34)]{OpenAI-2026c} in
\eqref{eq:kernel-compression}. The rapid decay of
$\widehat\psi_a$ gives the required logarithmic saving
for this choice of major arcs.

We next consider an edge from the indicator in
\eqref{eq:shift-kernel} that involves a transferred prime.
As in \cite[Section~4.6]{OpenAI-2026c}, we fix the new primes and choose
$p_*$ from the largest relevant set. With the real lattice
basis matrix $g\in\mathrm{SL}_2(\mathbb R)$ of
\cite[Section~3.4]{OpenAI-2026c}, each position has the form $g{\mathbf z}$,
where $ {\mathbf z}\in\mathbb Z^2$ is primitive. We choose
$w_{ {\mathbf z}}\in\mathbb Z^2$ with $\det( {\mathbf z},w_{ {\mathbf z}})=1$. We write the
final integer vector as $j_1w_{ {\mathbf z}}+b {\mathbf z}$ and put
$r=(gw_{ {\mathbf z}})_1/(g {\mathbf z})_1$. Modulo one, this is the root ratio
of $g {\mathbf z}$. The determinant condition is
\begin{equation}\label{eq:promoted-determinant}
 j_1=a(P_{\mathrm{full}}-Dp_*UT_f),
 \quad P_{\mathrm{full}}=DB.
\end{equation}
The distinctness conditions give $p_*\nmid P_{\mathrm{full}}$.
Since $p_*\nmid a$, we have
\begin{align*}
 j_1\equiv aP_{\mathrm{full}}\not\equiv0\ (\bmod\ p_*),
\end{align*}
independently of $D,U$. The line
$\ell_*\in\mathbb P^1(\mathbb F_{p_*})$ is attached to
the chosen entry in \cite[Section~4.6]{OpenAI-2026c}. The condition that $j_1w_{ {\mathbf z}}+b {\mathbf z}$ lie on this line modulo $p_*$
either has no solution or determines a unique residue class
$b\equiv b_*\;(\bmod p_*)$. The coordinate bounds give $|b+j_1r|\leq C_{\mathrm{box}}$,
where $C_{\mathrm{box}}$ depends only on $a$ and the fixed box.
We take the separation constant larger than $2C_{\mathrm{box}}$.
It follows from \eqref{eq:promoted-determinant} that
\begin{align*}
\Big\|aDUT_fr-\frac{aP_{\mathrm{full}}r+b_*}{p_*}
\Big\|_{ {\mathbb T}}
\leq \frac{C_{\mathrm{box}}}{p_*}.
\end{align*}
For newly selected primes satisfying the modified separation condition,
this inequality allows at most one value of $DU$, because
$p_*\geq\exp(L^{a_{i_*}})$ in the largest relevant set
$\mathcal B_{i_*}$. No prime dividing $U$ occurs in the source tuples.
The prime factors of $DU$ in those tuples are therefore exactly
the prime factors of $D$. Hence the product determines $D$ and
$U$ separately, and there is at most one choice of the omitted
entry in each tuple. Under the average over list positions,
such a tuple has weight $M_{\mathrm{pad}}^{-K}$, where
$M_{\mathrm{pad}}$ is the full list length in \cite{OpenAI-2026c}.

For this tuple, the bounds for the number of final positions and
the remaining multiplicity follow from
\cite[(4.23) and (4.37)]{OpenAI-2026c}.
The exceptional choices of new primes have exponentially small
total mass, with the same bound as in
\cite[Section~4.6]{OpenAI-2026c}. If we reverse an edge, both $T$ and
$B-A$ change sign, so the same bound holds for the adjoint
operator. For an edge involving a transferred prime in the
integral term of \eqref{eq:shift-kernel}, the kernel has
absolute value $O(L^{3A_0}/Y)$, where $A_0$ is the major-arc
exponent. The proof in \cite[Section~4.7]{OpenAI-2026c} depends only on
this bound, the lattice box and the congruence defining the final
line. The same estimate therefore holds for the modified kernel.

We apply these operator bounds in the moment argument of
\cite[Section~4.8]{OpenAI-2026c}. The distinctness conditions
on the initial prime labels remain the same. The indicator and
integral terms have $O_a(1)$ and $O_a(Y)$ possible final
vertices, respectively. The endpoint moment bounds and
the root-ratio estimates also allow us to remove the
admissibility restrictions, as in
\cite[(4.56)--(4.57)]{OpenAI-2026c}.

We next remove the auxiliary lists, keeping track of their normalization. Each list has
$J=M_{\mathrm{pad}}-1=\lfloor L^{0.01}\rfloor$ shared entries,
whose total product is $D$. On a block $d_0\leq D<2d_0$,
the two coordinate scales are $d_0YM$ and $N$, so their
product is $d_0YX$. The endpoint vectors are invariant under permutations of their
lists. Thus the permutation averages drop out of the expanded
pairing. After division by $d_0$, the kernel multiplier $d_0/D$
becomes $1/D$.
For the ordered shared labels $p_{i,1},\ldots,p_{i,J}$,
the normalization identity in \cite[(4.58)]{OpenAI-2026c} is
\begin{align*}
\frac1D\prod_{i=1}^K V_i^{-(J+2)}
=\Big(\prod_{i=1}^K\prod_{\ell=1}^J\frac1{V_ip_{i,\ell}}\Big)
 \prod_{i=1}^K V_i^{-2}.
\end{align*}
The sum of the first factor over unrestricted shared lists is one.
The second factor normalizes the two unshared prime products
in $Q_Y$. Each label occurs once in the ordered slots, so there
is no factorial or additional power of $M_{\mathrm{pad}}$.

For fixed disjoint unshared lists, the total harmonic weight
of shared lists that repeat a prime or meet an unshared
list is at most
\begin{align*}
O_K\Big((J+1)^2
\max_{\substack{1\leq i\leq K\\p\in\mathcal B_i}}\frac1{V_ip}\Big)
\ll_K L^{0.02}\exp(-L^{1/10}).
\end{align*}
We use absolute bounds to remove these restrictions from the signed
pairing. By \eqref{eq:affine-fourth-divisor}, the indicator contribution
is $O_{a,K}(XYL^{C_2})$. For the integral term,
the conditions
\begin{align*}
|Bv-Aw|\ll_aY,\quad A,B\asymp Y
\end{align*}
allow $O_a(1)$ values of either variable when the other
is fixed. The ordinary divisor moments and the kernel
bound $O(L^{3A_0}/Y)$ therefore give
$O_{a,K}(XYL^{C_2+3A_0})$, after increasing $C_2$ if necessary.
The exponent $C_2$ depends on the fixed coefficient bounds
and is independent of $K$. By the preceding bound on the harmonic
weights, the error from repeated or forbidden shared labels is
smaller than $XY$ times every fixed negative power of $L$.
The same bound for the absolute integral allows us to restore
pairs of unshared products with a common prime. Their normalized
count is $O_K(Y^2\exp(-L^{1/10}))$.

We also keep track of the scale in the signed estimate.
With accuracy $E_1$ in the moment bound, the endpoint estimates give
$O_{a,K}(d_0YX L^{-E_1+C_*})$ on the block
$d_0\leq D<2d_0$, where $C_*$ depends only on the fixed
coefficient bounds and is independent of $A_0$ and $K$.
After division by $d_0$, this becomes
$O_{a,K}(XYL^{-E_1+C_*})$. There are $O_K(L)$ blocks,
and every ordered shared list belongs to exactly one of
them. Their sum is therefore
\begin{align*}
O_{a,K}(XYL^{-E_1+C_*+1}).
\end{align*}
For a prescribed replacement accuracy $E_0>0$, we first choose
$E_1>E_0+C_*+2$. We then fix the major-arc exponent $A_0$
and finally $K$. The sum over blocks costs one power of $L$;
our choice of accuracy absorbs this loss before we choose $K$.
It follows that $Q_Y$ differs from its major-arc contribution
by $O(XYL^{-E_0})$.

Finally, on each major arc, we write
$\theta=b/q_0+\lambda/Y$, as in \cite[Section~5]{OpenAI-2026c}.
The additional phase factors as
\begin{align*}
e\Big(a\theta(A-B)\Big)=e(a\theta A)e(-a\theta B).
\end{align*}
The rational phase in the six variables has a character expansion
on reduced residue classes. The sum of the absolute values of
its coefficients is at most $q_0^6$. The smooth product weights are $\eta(v)e(\pm a\lambda v)$.
Their derivative bounds involve the same powers of $L$, since
$a$ is fixed.
We use $\psi_a$ in the Mellin kernel and retain the Dirichlet
polynomials formed from $\alpha_m$ and $\beta_n$ in the
original proof. The estimates in
\cite[Lemmas~5.2--5.3 and Proposition~5.1]{OpenAI-2026c} then give
a major-arc contribution of $O(XYL^{-E_0})$.
We choose $E_0>2A+4$ and combine this bound with
$|S_Y|^2\ll(X/Y)Q_Y$. There are $O_K(L)$ dyadic values of $Y$,
so their sum satisfies
\eqref{eq:type-ii}.

The dependence on $a$ and the fixed sets $\mathcal B_i$ enters
only through the implied constants and the threshold for $x$.
Hence the required lower bound for $K$ depends only on the
fixed parameters in the lemma.
\end{proof}

\subsection{A correlation estimate with an invariant factor}
For the centered correlation estimate, we choose two families
of pairwise disjoint prime sets $\mathcal P_g$.
Each set consists of all primes in an interval, with
$1/2\leq\sum_{p\in\mathcal P_g}p^{-1}\leq2$.
We require
\begin{equation}\label{eq:bands}
\begin{aligned}
\mathcal P_g
&\subset\{p:\exp(L^{0.245})<p<\exp(L^{0.249})\}
&&\text{in the first family},\\
\mathcal P_g
&\subset\{p:\exp(L^{0.34})<p<\exp(L^{0.39})\}
&&\text{in the second family}.
\end{aligned}
\end{equation}
For example, we may take the sets
$\{p:e^j<\log p\leq e^{j+1}\}$ for integers $j$ whose
corresponding prime intervals lie within the indicated ranges.
Mertens' theorem gives a reciprocal prime sum of $1+o(1)$
for each set, and this construction yields $\asymp\log L$
sets in each family. These ranges satisfy all the size and
reciprocal-mass conditions in \cite[(2.1)]{Walsh-2026}.
Their union $\mathcal C$ satisfies
$\mathcal C\subseteq\mathcal S$ and
$\mathcal B\cap\mathcal C=\varnothing$, and every prime in
$\mathcal C$ exceeds $Z$ for sufficiently large $x$. For each $\mathcal P_g$, we define
\begin{align*}
V_g=\sum_{p\in\mathcal P_g}p^{-1},
\quad
\omega_g(n)=\sum_{p\in\mathcal P_g}\mathbf 1_{p\mid n}.
\end{align*}
For every set $\mathcal P\subseteq\mathcal S$, we define
\begin{align*}
\alpha_{\mathcal P}(n)=\tau_k(n_{\mathcal P'})
-\tau_k(n_{\mathcal S\setminus\mathcal P})
  \frac{M_k(n_{\mathcal S'})}{V}.
\end{align*}
As in \cite[(2.10)]{Walsh-2026}, the definition gives
$\alpha_{\mathcal P}(pn)=\alpha_{\mathcal P}(n)$ for $p\in\mathcal P$.
Multiplicativity yields
\begin{equation}\label{eq:invariance}
\xi_k(n)=\tau_k(n_{\mathcal P})\alpha_{\mathcal P}(n).
\end{equation}
For a selected subfamily indexed by $\mathcal G$ and
$\mathcal P=\bigcup_{g\in\mathcal G}\mathcal P_g$, we also put
\begin{align*}
\Omega_{\mathcal P}(n)=\prod_{g\in\mathcal G}
 V_g^{-1}\mathbf 1_{\omega_g(n)=1}.
\end{align*}

The following proposition applies to an invariant second
factor bounded by a divisor function, as required for the
convolution of sifted factors in Section~\ref{sec:composites}.

\begin{proposition}\label{prop:graph}
The positive integer $j$ and the constants $c,C>0$ are fixed.
The selected subfamily contains at least $c\log L$ sets from
each family in \eqref{eq:bands}. The weight $H$ satisfies
$H(n)=H(n_{\mathcal S\setminus\mathcal P})$ and $|H(n)|\leq L^C$,
and the second factor $B$ satisfies
\begin{align*}
 B(pn)=B(n)\quad(p\in\mathcal P),
 \quad |B(n)|\le L^C \tau_j(n_{\mathcal P'}).
\end{align*}
For every fixed $A>0$, we have
\begin{equation}\label{eq:graph}
 \sum_n\Phi(n/x)H(n)\alpha_{\mathcal P}(n)B(n-h)
 \Omega_{\mathcal P}(n)\Omega_{\mathcal P}(n-h)
 \ll_A xL^{-A},
\end{equation}
uniformly in the selected subfamilies and in the functions satisfying these conditions.
\end{proposition}

\begin{proof}
We use \cite[Theorem~7.2]{Walsh-2026} together with the analytic
estimates in the proof of \cite[Proposition~4.1]{Walsh-2026}.
We set $A_0(n)=H(n)\alpha_{\mathcal P}(n)$.
We check the moment bounds and invariance of the endpoint weights
$A_0$ and $B$, and then prove the required cancellation.

First, we fix integers $r,m_0\geq1$ and a scale $Y=x^{1+o(1)}$.
The collection $\mathbf p$ consists of one ordered list of $m_0$
distinct primes from each selected set. The product of all
entries is denoted by $P_0$.
Our prime ranges satisfy the hypotheses of
\cite[Lemma~2.3]{Walsh-2026}. The bounds on $H$ and $B$ cost an additional factor $L^{rC}$
in the moment estimate. The proof of that lemma therefore gives
\begin{equation}\label{eq:graph-moments}
\begin{aligned}
\sum_{\mathbf p}\Big(\prod_{g\in\mathcal G}V_g^{-m_0}\Big)
\sum_{\substack{n\asymp Y\\P_0\mid n}}
\Big(|A_0(n)|^r+|B(n)|^r\Big)
\ll YL^{r(C+k-1)+\max(k,j)^r-1}.
\end{aligned}
\end{equation}
The exponent of $L$ and the implied constant do not depend on
$m_0$, although the threshold for $x$ may depend on $m_0$.
The cases $r=2,4$ give the moment bounds needed below.

Second, both endpoint weights are invariant under multiplication
by primes in $\mathcal P$ and under permutations within each
selected list. Together with \eqref{eq:graph-moments}, these
properties allow us to apply \cite[Theorem~7.2]{Walsh-2026}.
We pass from the graph pairing to the arithmetic comparison
by removing the shared prime labels. The argument below also
bounds the error from repeated labels.

For the comparison, we fix nonnegative integers $t_g$.
The collection $\mathbf p$ now consists of ordered lists.
Its $g$th list contains $t_g$ distinct primes from $\mathcal P_g$.
For a function $b$ of these lists, we define
\begin{align*}
W_{\mathbf t}^b(n)=
\Big(\prod_{g\in\mathcal G}
\mathbf 1_{\omega_g(n)=t_g}V_g^{-t_g}\Big)
\sum_{\substack{\mathbf p\\
p_{g,i}\mid n\ \text{for all }g,i}}b(\mathbf p).
\end{align*}
In particular, $W_{\mathbf1}^{1}=\Omega_{\mathcal P}$.
We set $\Psi(t)=t\Phi(t)$ and write
\begin{align*}
C_0=\sum_n\frac{\Psi(n/x)}{n}
A_0(n)B(n-h)
\Omega_{\mathcal P}(n)\Omega_{\mathcal P}(n-h).
\end{align*}
The construction in \cite[Proposition~6.1]{Walsh-2026}
uses the moments and invariance of the endpoint weights.
Its analytic estimates require cancellation only at the first
endpoint. No centering condition on $B$ is needed.
For every fixed $D>0$, the construction gives a fixed integer
$m_0\geq1$ and finitely many comparison terms $C_\nu$ such that
\begin{equation}\label{eq:graph-comparison}
C_0+\sum_\nu c_\nu C_\nu=O(L^{-D}),
\quad
\#\{\nu\}+\sum_\nu|c_\nu|\leq L^{C_1},
\end{equation}
where $C_1$ is fixed once the graph parameters have been chosen.

We describe these terms as follows.
The second family is partitioned into two parts, each containing
$\gg\log L$ selected sets.
For each comparison $\nu$, the integers
$r_{\nu,g}\in\{0,\ldots,m_0\}$ satisfy $r_{\nu,g}=0$
in the first family. In each part of the second family,
at least one $r_{\nu,g}$ is positive.
For each $g$, the sum runs over ordered tuples
$(d_{g,1},\ldots,d_{g,r_{\nu,g}})\in\mathcal P_g^{r_{\nu,g}}$.
The entries in these auxiliary tuples need not be distinct.
We denote the collections from the two parts by
$\mathbf d_1,\mathbf d_2$ and their respective products by
$D_1,D_2$. We set $t_{\nu,g}=1+r_{\nu,g}$.

The removal of the shared labels follows
\cite[Section~7.5]{Walsh-2026}, with shift $-h$.
We denote their product by $D_R$. The original graph
endpoints satisfy $n'=n-hD_RD_1D_2$, and the change
of variables is
\begin{align*}
n=D_Rw,\quad n'=D_R(w-hD_1D_2).
\end{align*}
The endpoint invariance removes $D_R$ from $A_0(n)$
and $B(n')$, while the harmonic factor satisfies
\begin{align*}
\frac{\Psi(n/(xD_RD_1D_2))}{n}
=\frac1{D_R}\frac{\Psi(w/(xD_1D_2))}{w}.
\end{align*}
Here $M_{\mathrm{gr}}$ is the full-list length in the graph
construction, fixed sufficiently large. The product of shared
labels in the $g$th set is $(D_R)_g$. The local normalization is
\begin{align*}
\frac{V_g^{-M_{\mathrm{gr}}-t_{\nu,g}}}{(D_R)_g}
=V_g^{-2t_{\nu,g}}
 \prod_{\substack{p\text{ in the shared}\\
                  \text{slots of }\mathcal P_g}}\frac1{V_gp}.
\end{align*}
Thus the shared labels have normalized harmonic weights.
The remaining factors give the two complete-list weights.

For now, the auxiliary primes are still required to avoid the
unshared endpoint lists.
For fixed reduced endpoints and auxiliary lists, there are
$O_{M_{\mathrm{gr}}}(\log L)$ shared slots. Each prime
has harmonic weight at most $2\exp(-L^{0.245})$.
On the support of the smooth weight, there are
$O_{M_{\mathrm{gr}}}(L)$ primes that divide
$w(w-hD_1D_2)D_1D_2$ or occur in another shared slot.
Hence the total weight of shared lists with a repeated label
or one of these divisibilities is
\begin{align*}
\ll L^{O_{M_{\mathrm{gr}}}(1)}\exp(-L^{0.245}).
\end{align*}
If each reduced endpoint has exactly $t_{\nu,g}$ prime
divisors in every active set, all remaining choices of
shared labels reconstruct admissible full lists.
If a reduced endpoint has more than $t_{\nu,g}$ divisors
in some set, an admissible full list must contain a shared
prime dividing that endpoint. Such choices are covered
by the same bound. The total weight of valid shared lists
is consequently
\begin{align*}
\prod_{g\in\mathcal G}
\mathbf 1_{\omega_g(w)=\omega_g(w-hD_1D_2)=t_{\nu,g}}
+O\Big(L^{O_{M_{\mathrm{gr}}}(1)}\exp(-L^{0.245})\Big).
\end{align*}
This statement applies on the original graph range, where
both reduced counts lie between $t_{\nu,g}$ and
$M_{\mathrm{gr}}$.

The counts in each set are bounded by the fixed number
$M_{\mathrm{gr}}$. There are $O(\log L)$ sets, so the absolute
sum over the unshared lists costs at most a fixed power of $L$.
At scale $Y=xD_1D_2$, the standard second moment estimates give
\begin{align*}
\frac1Y\sum_{w\asymp Y}
|A_0(w)B(w-hD_1D_2)|\ll L^{O(1)}.
\end{align*}
The error from the shared lists is therefore smaller than every
fixed negative power of $L$. This remains true after summation
against the comparison coefficients.
The fourth moment estimate below gives the same conclusion
for the error from removing the auxiliary-prime exclusions.
After both errors are absorbed, the comparison term is
\begin{equation}\label{eq:graph-comparison-term}
\begin{aligned}
C_\nu={}&\sum_{\mathbf d_1,\mathbf d_2}
a_\nu(\mathbf d_1)b_\nu(\mathbf d_2)
\Big(\prod_{g\in\mathcal G}\prod_{i=1}^{r_{\nu,g}}
\frac1{V_gd_{g,i}}\Big)
\sum_w\frac{\Psi(w/(xD_1D_2))}{w}\\
&\times A_0(w)B(w-hD_1D_2)
W_{\mathbf t_\nu}^{u_\nu}(w)
W_{\mathbf t_\nu}^{v_\nu}(w-hD_1D_2).
\end{aligned}
\end{equation}
The product weights sum to one because
$\sum_{p\in\mathcal P_g}(V_gp)^{-1}=1$.
This comparison has the form of \cite[(6.4)]{Walsh-2026},
with shift $-h$. The coefficients satisfy
$|a_\nu|,|b_\nu|,|u_\nu|,|v_\nu|\leq1$.
The functions $u_\nu$ and $v_\nu$ depend only on the lists
from the second family. These are the bounds and restrictions
in \cite[Proposition~6.1]{Walsh-2026}.

We remove the exclusions at auxiliary primes by the argument of
\cite[Section~7.5]{Walsh-2026}. The auxiliary lists are fixed,
and $Y=xD_1D_2$.
For every prime $p\mid D_1D_2$, the conditions $p\mid w$
and $p\mid w-hD_1D_2$ are equivalent.
Both arguments are comparable to $Y$ on the smooth support,
since $hD_1D_2/Y=h/x$.
The pointwise bounds and the standard fourth moment estimates
give, uniformly in the auxiliary lists,
\begin{align*}
\sum_{w\asymp Y}|A_0(w)|^4
\ll YL^{C_4},\;
\sum_{w\asymp Y}|B(w-hD_1D_2)|^4\ll YL^{C_4},
\end{align*}
where $C_4$ is fixed.
For the second estimate, we translate the summation variable
and apply the usual moment bound. There are $O(Y/p)$ multiples
of $p$ in the summation range, since $Y/p\to\infty$ uniformly.
H\"older's inequality now gives
\begin{align*}
&\frac1Y\sum_{\substack{w\asymp Y\\p\mid w}}
|A_0(w)B(w-hD_1D_2)|\\
&\ll \frac1Y\Big(\frac Yp\Big)^{1/2}
\Big(\sum_{w\asymp Y}|A_0(w)|^4\Big)^{1/4}
\Big(\sum_{w\asymp Y}|B(w-hD_1D_2)|^4\Big)^{1/4}\ll L^{C_4/2}p^{-1/2}.
\end{align*}
The smooth harmonic factor $\Psi(w/Y)/w$ is $O(Y^{-1})$.
On the support of each list weight, the $t_{\nu,g}$ primes
from $\mathcal P_g$ admit $t_{\nu,g}!$ orderings.
The absolute value of the product of the two list weights
is therefore at most
\begin{align*}
\prod_{g\in\mathcal G}(t_{\nu,g}!)^2V_g^{-2t_{\nu,g}}
=L^{O_{ {m_0}}(1)}.
\end{align*}
There are $O_{ {m_0}}(\log L)$ auxiliary labels, all corresponding
to primes exceeding $\exp(L^{0.34})$. We sum over these
labels and obtain a total collision error of
$L^{O_{ {m_0}}(1)}\exp(-\tfrac12L^{0.34})$.
This is smaller than every fixed negative power of $L$,
even after we average over the auxiliary lists and sum
against the comparison coefficients. Finally, any prime common
to the two complete endpoint lists divides $hD_1D_2$. For sufficiently large $x$, such a prime
exceeds $|h|$ and must therefore divide $D_1D_2$.
The preceding estimate already covers this case.

Third, we prove cancellation in the factor $A_0$ in each
comparison. We restrict $D_i$ to a dyadic interval
$H_i\leq D_i<2H_i$ and use the scale
$Y=xH_1H_2\asymp xD_1D_2=x^{1+o(1)}$, as in
\cite[Section~6.1]{Walsh-2026}. The estimates in \cite[Section~6.1]{Walsh-2026} require only
the $L^2$ norm of $B$. The case $r=2$ of
\eqref{eq:graph-moments} gives this bound. It suffices for the
bilinear Fourier estimate and the local Mellin reduction of
\eqref{eq:graph-comparison-term}.

For low Mellin frequencies, we fix 
$A=n_{\mathcal S}$ and write $n=At$. The definition of
$\alpha_{\mathcal P}$ gives
\begin{align*}
\alpha_{\mathcal P}(n)
=\tau_k(A_{\mathcal S\setminus\mathcal P})
\Big\{\tau_k(t)-\frac{M_k(t)}{V}\Big\},
\quad (t,\mathcal S)=1.
\end{align*}
Moreover, $H(n)=H(A_{\mathcal S\setminus\mathcal P})$,
and the list weights depend only on $A$. The coefficient of the centered inner sum has absolute value
at most $L^{C_2}\tau_k(A)$. The periodic modulus in the Mellin
reduction is at most a fixed power of $L$. It is therefore
$\mathcal S$-smooth for sufficiently large $x$. The smooth
cutoffs have logarithmic derivatives bounded by fixed powers of
$L$.  {We keep the original small-prime parameters fixed and apply
\cite[Lemma~5.2]{Walsh-2026} at scale $Y=x^{1+o(1)}$.
The lemma bounds the sum over $A$ of the absolute values of
the centered sums. After normalization, this bound is}
$L^{C_3}\exp(-c_3L^{1/10})$, uniformly for Mellin
frequencies bounded by any fixed power of $L$.

For high Mellin frequencies, we fix a set $\mathcal P_g$
from the first family in \eqref{eq:bands}.
The corresponding list has length one in every comparison.
The comparison coefficient is independent of the chosen prime.
On the support of the list weight, $n$ has exactly one prime
divisor $p$ in $\mathcal P_g$.

Every $p\in\mathcal P_g$ satisfies $p^2<Y$.
Hence the number of integers $n\asymp Y$ divisible by the square
of a prime in $\mathcal P_g$ is
$O(Y\sum_{p\in\mathcal P_g}p^{-2})$.
Cauchy--Schwarz and the fourth moment estimate give
\begin{align*}
\sum_{\substack{n\asymp Y\\
p^2\mid n\ \text{for some }p\in\mathcal P_g}}
|A_0(n)|^2\ll YL^{O(1)}
\Big(\sum_{p\in\mathcal P_g}p^{-2}\Big)^{1/2}
\ll YL^{O(1)}\exp(-cL^{0.245}).
\end{align*}
The list weights contribute at most a fixed power of $L$.
The mean value estimate for Dirichlet polynomials used below
shows that this contribution is negligible after Mellin integration.
We omit these terms and write $n=pm$ for the remaining integers.
The invariance identity
\begin{align*}
H(pm)\alpha_{\mathcal P}(pm)
=H(m)\alpha_{\mathcal P}(m)
\end{align*}
then gives the prime Dirichlet polynomial in
\cite[Section~6.4]{Walsh-2026}.
A dyadic decomposition of the prime range and character
orthogonality give polynomials of the form
\begin{equation}\label{eq:prime-polynomial}
\frac1P\sum_{p\in\mathcal P_g}f(p/P)\chi(p)p^{it},
\end{equation}
where $f$ is a fixed smooth dyadic cutoff supported in $[1/2,3]$.
The modulus of $\chi$ is bounded by a fixed power of $L$.

We first check the length condition for the remaining
Dirichlet polynomial.
At scale $Y\asymp xD_1D_2$, its length is $Y/P$.
The frequencies satisfy $|t|\leq T_0=xL^{B_1}$,
where $B_1>0$ is fixed.
Each of $D_1$ and $D_2$ contains a prime from the second family.
Hence \eqref{eq:bands} gives
\begin{equation}\label{eq:graph-length}
\frac{T_0}{Y/P}
\ll \frac{PL^{B_1}}{D_1D_2}
\ll L^{B_1}\exp\{L^{0.249}-2L^{0.34}+O(1)\}
=o(1).
\end{equation}

We next smooth the restriction $p\in\mathcal P_g$
in \eqref{eq:prime-polynomial}, as in
\cite[Section~6.4]{Walsh-2026}.
We set $\varepsilon=L^{-K_0}$, where the fixed integer $K_0$
will be chosen below.
We choose a smooth cutoff $w_P$ with $0\leq w_P\leq1$.
The factor $w_P(p/P)$ vanishes outside the interval defining
$\mathcal P_g$. Inside this interval, it equals one except
on intervals of total length $O(\varepsilon P)$ near the endpoints.
For every fixed $r$, we require
$\|w_P^{(r)}\|_\infty\ll_r\varepsilon^{-r}$,
with derivatives taken in the scaled variable.

We denote \eqref{eq:prime-polynomial} by $Q(t)$.
The smoothed polynomial $Q_{\mathrm{sm}}(t)$ has $fw_P$
in place of $f$.
The intervals near the endpoints contain $O(\varepsilon P+1)$
integers. Thus
\begin{align*}
|Q(t)-Q_{\mathrm{sm}}(t)|
\ll \varepsilon+P^{-1}
\end{align*}
uniformly in $t$.

We use this pointwise bound together with the mean square
estimate for the remaining polynomial. After the separation of
variables in \cite[Section~6.4]{Walsh-2026}, the latter has the form
\begin{align*}
R(t)=\frac PY\sum_{m\asymp Y/P}a_m m^{it},
\quad
\sum_{m\asymp Y/P}|a_m|^2\ll (Y/P)L^{C_5}.
\end{align*}
For each fixed choice of the separated parameters,
the coefficients $a_m$ are independent of $t$.
The bound $|A_0(m)|\ll L^{C+k-1}\tau_k(m)$ and the fixed
logarithmic bound for the remaining list weight reduce the
coefficient estimate to the usual second moment of $\tau_k$.
The estimate remains valid after restriction to residue classes
or multiplication by factors of modulus one.

By \eqref{eq:graph-length} and the mean value estimate for
Dirichlet polynomials, we have
\begin{align*}
&\int_{-T_0}^{T_0}|R(t)|^2\,dt
\ll
\frac{T_0+Y/P}{(Y/P)^2}
\sum_{m\asymp Y/P}|a_m|^2
\ll L^{C_5},\\
&\int_{-T_0}^{T_0}
\Big|(Q(t)-Q_{\mathrm{sm}}(t))R(t)\Big|^2\,dt
\ll L^{C_5}(\varepsilon+P^{-1})^2.
\end{align*}
The estimates are uniform in the separated Mellin frequencies.
Minkowski's inequality bounds the integrals against the smooth
kernels. The subsequent sums over dyadic intervals and residue
classes cost at most a fixed power of $L$.
After these losses and the desired saving are fixed, we choose
$K_0$ sufficiently large. Since $P\geq\exp(L^{0.245})/2$,
the smoothing error has the required logarithmic saving.

For this choice of $K_0$, every fixed-order derivative of $fw_P$
is bounded by a fixed power of $L$. The support lies inside
the interval defining $\mathcal P_g$, so the restriction
$p\in\mathcal P_g$ may be omitted.
For sufficiently large $x$, every prime in the support exceeds
the character modulus.
We apply the full estimate of
\cite[Lemma~5.1]{Walsh-2026}.
For sufficiently large $x$, we have
$P\geq\exp(L^{0.245})/2\geq\exp(L^{0.24})$, which verifies
the length condition. The height condition covers the full
range $|t|\leq T_0=xL^{B_1}<x^2$.
For every fixed integer $\ell\geq1$, the lemma gives
\begin{align*}
|Q_{\mathrm{sm}}(t)|
\ll_\ell L^{O_\ell(1)}(1+|t|)^{-\ell}
+P^{-1/16}L^{O(1)}.
\end{align*}
The second term is smaller than every fixed negative power of $L$.
The rapid decay of the smooth kernels allows us first to restrict
the additional Mellin frequencies to fixed logarithmic ranges.
We then choose the cutoff between low and high frequencies.
This cutoff is a fixed power of $L$.
The power is chosen large enough to dominate the additional
frequency ranges and absorb the derivative bounds in the
preceding estimate.
Above this cutoff, the shifted frequencies have absolute value
comparable to $|t|$.
The bound therefore gives the required decay after separation
of variables.
We combine this decay with the mean square bound for $R(t)$
to obtain the required logarithmic saving.

We take $D\geq A+1$ in \eqref{eq:graph-comparison}.
Once $C_1$ is fixed, we choose enough logarithmic saving
in the analytic estimates to absorb the factor $L^{C_1}$.
The sum of the comparison bounds then gives $C_0\ll L^{-A}$.
Finally, the identity $\Psi(n/x)/n=\Phi(n/x)/x$ gives
\eqref{eq:graph} after multiplication by $x$.
\end{proof}

\section{The sum over composite integers}\label{sec:composites}

We first fix $s$ and $\gamma$ as in Proposition~\ref{prop:presieve}.
We then fix the truncation heights and the required accuracy
in the Type II estimate. The number $K$ and the exponents
defining the prime sets of $\mathcal W$ are chosen afterwards.
All auxiliary parameters remain fixed as $x\to\infty$.

For sufficiently large $x$, a composite integer $v\asymp x$
with $P^-(v)>x^\gamma$ has between two and $2/\gamma$ prime
factors, counted with multiplicity. We shall prove that,
for each integer $2\leq j\leq2/\gamma$,
\begin{equation}\label{eq:composites}
\sum_{p_1,\ldots,p_j>x^\gamma}
\Phi((p_1\cdots p_j+h)/x)
\xi_k(p_1\cdots p_j+h)
\mathcal W(p_1\cdots p_j+h)
=o(xL^{k-2}).
\end{equation}
Throughout this section, $p_1,\ldots,p_j$ denote primes.

We first discard the tuples with repeated primes. Each
corresponding product is divisible by $p^2$ for some
$p>x^\gamma$, so there are only $O(x^{1-\gamma})$
possible products. For fixed $j$, the factorization
multiplicities and divisor weights together contribute
a factor $O_\varepsilon(x^\varepsilon)$. The total
contribution is therefore $O_\varepsilon(x^{1-\gamma+\varepsilon})$.
We choose $0<\varepsilon<\gamma$ to obtain an error of
$o(xL^{-A})$ for every fixed $A>0$.

\subsection{Replacement of the prime indicators}
We define $a(1)=0$ and
\begin{equation}\label{eq:prime-replacement-weight}
a(m)=
\frac{\mathbf 1_{m>x^\gamma}\mathbf 1_{P^-(m)>Z}}
{V(Z)\log m},
\quad m\geq2.
\end{equation}
 {We take $\mathcal P=\mathcal B$ in the definition of}
$\alpha_{\mathcal P}$. If $p^2\nmid n$ for every
$p\in\mathcal B$, the identity \eqref{eq:invariance} gives
\begin{equation}\label{eq:damped-factorization}
\xi_k(n)\mathcal W_\rho(n)
=k^K\alpha_{\mathcal B}(n)\mathcal W_{k\rho}(n),
\quad k\rho=\frac12.
\end{equation}
We combine the estimate
\begin{align*}
\sum_{p\in\mathcal B}p^{-2}
\ll_K\exp(-L^{1/10})
\end{align*}
with the fixed divisor moment bounds to discard the terms
for which $p^2\mid n$ for some $p\in\mathcal B$.
Their contribution is $O(xL^{-A})$ for every fixed $A>0$.

To obtain a bounded weight for Lemma~\ref{lem:type-ii},
we define
\begin{align*}
F(n)=L^{-B_0}\alpha_{\mathcal B}(n)
\mathbf 1_{|\alpha_{\mathcal B}(n)|\leq L^{B_0}}.
\end{align*}
Then $|F(n)|\leq1$, and the invariance of
$\alpha_{\mathcal B}$ gives $F(pn)=F(n)$ for every
$p\in\mathcal B$. We replace $\alpha_{\mathcal B}$
by $L^{B_0}F$ and estimate the resulting error using
\begin{align*}
|\alpha_{\mathcal B}(n)|
\mathbf 1_{|\alpha_{\mathcal B}(n)|>L^{B_0}}
\leq L^{-B_0}|\alpha_{\mathcal B}(n)|^2.
\end{align*}
We bound the factorization multiplicities by $\tau_j(v)$
and apply the positive majorant and the fixed divisor
moment estimates to obtain
\begin{equation}\label{eq:clipping-tail}
L^{-B_0}\sum_{v\asymp x}
Q_k(v+h)^2\tau_j(v)L^{C_j}
\ll xL^{C_{k,j,2}-B_0}.
\end{equation}
The exponent $C_{k,j,2}$ is independent of $K$, although
the implied constant may depend on the fixed sets
$\mathcal B_i$. We may therefore choose
$B_0>C_{k,j,2}+A$ before choosing $K$, so that
\eqref{eq:clipping-tail} is $O(xL^{-A})$.
A weighted second-moment bound suffices for this truncation.

We now replace the $j$ prime indicators, one at a time, by
the coefficients in \eqref{eq:prime-replacement-weight}.
In each difference, let $m$ denote the factor being replaced
and $n$ the product of the remaining $j-1$ factors.
We decompose their ranges into dyadic intervals, with
$M,N\geq x^{\gamma/2}$ and $MN\asymp x$.
The coefficient of $n$ vanishes unless $P^-(n)>Z$,
and its absolute value is bounded by $\tau_{j-1}(n)$ times
a fixed power of $L$. We truncate this coefficient at
a sufficiently large power of $L$. The earlier truncation
of $\alpha_{\mathcal B}$ and the fixed moments of
$\tau_{j-1}$ ensure that the resulting error is negligible.

To separate $m$ and $n$ in the smooth weight, we define
\begin{align*}
\widehat\Phi_{h,x}(t)=\frac1{2\pi}
\int_0^\infty\Phi(v+h/x)v^{-it}\,\frac{dv}{v}.
\end{align*}
Mellin inversion gives
\begin{equation}\label{eq:shifted-mellin}
\Phi((mn+h)/x)
=\int_{\mathbb R}\widehat\Phi_{h,x}(t)
x^{-it}m^{it}n^{it}\,dt.
\end{equation}
For every fixed $B>0$, we have
$\widehat\Phi_{h,x}(t)\ll_{B,\Phi}(1+|t|)^{-B}$
uniformly for sufficiently large $x$.
We insert \eqref{eq:shifted-mellin} and combine this decay
with the same moment bounds to truncate the integral at
a fixed power of $L$. The condition $m>x^\gamma$
restricts $m$ to a subinterval of its dyadic interval,
as permitted in Lemma~\ref{lem:type-ii}.

We apply that lemma with $q=k\rho=1/2$. We choose the truncation
heights and the required logarithmic accuracy first, and then choose
$K$ sufficiently large for every $2\leq j\leq2/\gamma$.
By \eqref{eq:damped-factorization}, we obtain
\begin{equation}\label{eq:replacement}
\begin{aligned}
 &\sum_{p_1,\ldots,p_j>x^\gamma}
  \Phi((p_1\cdots p_j+h)/x)\xi_k(p_1\cdots p_j+h)
  \mathcal W(p_1\cdots p_j+h)\\
 &=\sum_{m_1,\ldots,m_j}
  \Phi((m_1\cdots m_j+h)/x)\xi_k(m_1\cdots m_j+h)
  \mathcal W(m_1\cdots m_j+h)\prod_{i=1}^j a(m_i)
  +o(xL^{k-2}).
\end{aligned}
\end{equation}
The factor $k^K$ and the bound for $\mathcal W$ affect the implied
constants in \eqref{eq:replacement}; they do not affect the
logarithmic exponents used to choose $K$.

\subsection{Invariant coefficients on sifted integers}
 {With $\mathcal C$ as defined after \eqref{eq:bands},}
for $m\geq1$ we write
$\bar m=m_{\mathcal C'}$ and define
\begin{equation}\label{eq:rough-factor}
\widetilde a(m)=
\begin{cases}
\displaystyle
\frac{\mathbf 1_{P^-(m)>Z}}{V(Z)\log\bar m},
&\bar m>x^\gamma,\\[5pt]
0,&\bar m\leq x^\gamma.
\end{cases}
\end{equation}
Since every prime in $\mathcal C$ exceeds $Z$, the definition
gives $\widetilde a(pm)=\widetilde a(m)$ for all
$p\in\mathcal C$. We shall use the following uniform upper
bound to show that we may replace $a$ by $\widetilde a$
in \eqref{eq:replacement} with a negligible error.

\begin{lemma}\label{lem:affine}
The positive integers $a,b$ and the positive constants
$\delta,C$ are fixed. The scale $R$ and the positive integer
$N$ satisfy $R\geq x^\delta$, $N\leq x^C$ and $P^-(N)>Z$.
For sets of primes
$\mathcal U_i,\mathcal V_i$, we define
\begin{align*}
 f_1(n)=\tau_a(n_{\mathcal V_1})\mathbf 1_{(n,\mathcal U_1)=1},\; f_2(n)=\tau_b(n_{\mathcal V_2})\mathbf 1_{(n,\mathcal U_2)=1}.
\end{align*}
Then, for sufficiently large $x$,
\begin{equation}\label{eq:affine-rough}
 \sum_{R<m\leq2R}f_1(m)f_2(Nm+h)
 \ll_{a,b,h,\delta,C}R\prod_{p\leq R}
 \int_{\mathbb Z_p}f_1(p^{v_p(t)})f_2(p^{v_p(Nt+h)})\,dt.
\end{equation}
For coprime integers $A,B\leq E$ with $(AB,h)=1$
whose prime factors all exceed $\exp(L^{1/10})$,
the unique integer $r$ satisfying
$-AB<r\leq0$, $A\mid r$ and $B\mid r-h$ determines the forms
\begin{equation}\label{eq:record-forms}
 \ell_1(t)=Bt+r/A,\quad \ell_2(t)=At+(r-h)/B.
\end{equation}
Then
\begin{equation}\label{eq:affine-product}
\begin{aligned}
 \sum_{x<ABm+r\leq2x}f_1(\ell_1(m))f_2(\ell_2(m))\ll_{a,b,h}\frac{x}{AB}\prod_{p\leq x/(AB)}
 \int_{\mathbb Z_p}
 f_1(p^{v_p(\ell_1(t))})f_2(p^{v_p(\ell_2(t))})\,dt,
\end{aligned}
\end{equation}
uniformly in $A,B$ and the sets of primes.
\end{lemma}

\begin{proof}
We use \cite[New Theorem~5, p.~377]{Henriot-2014}, which corrects
\cite[Theorem~5]{Henriot-2012}. For \eqref{eq:affine-rough},
the linear forms are $T$ and $NT+h$.
Since $P^-(N)>Z$ and $h$ is fixed, we have $(N,h)=1$
for sufficiently large $x$. Hence both forms are primitive. Their product
$Q(T)=T(NT+h)$ has degree two and height $O_h(N+1)$.
We use interval exponent $1/2$ and height exponent
$0<\eta<\min(1,\delta/(2C+2))$ in the cited theorem.
Since $R\geq x^\delta$ and $N\leq x^C$, we have
$R/\|Q\|^\eta\to\infty$ uniformly in $N$.

The functions $f_i$ are nonnegative and submultiplicative.
They satisfy
\begin{align*}
f_i(p^\nu)\leq\max(a,b)^\nu,
\quad f_i(n)\ll_{a,b,\varepsilon}n^\varepsilon,
\end{align*}
uniformly in the prime sets, for every fixed $\varepsilon>0$.
We choose $\varepsilon$ small enough to satisfy the growth
condition in the corrected theorem.

We denote by $\rho(p)$ the number of roots of $Q(T)$ modulo $p$,
and by $I_p$ the local integral in \eqref{eq:affine-rough}.
Since the forms are primitive, $\rho(p)\leq2$.
The local integral satisfies $I_p\geq1-\rho(p)/p$.
We use this inequality to compare the local factors in the
corrected theorem with $I_p$.

The density in \cite[(0.1)--(0.2)]{Henriot-2014} specifies
the exact valuations of both forms whenever at least one is
positive. The term with both exponents zero imposes no
congruence condition and is equal to $1$.
All terms are nonnegative, so we may sum over all exponent
choices at primes up to the theorem's cutoff.
The resulting local factor is
\begin{align*}
1+I_p-(1-\rho(p)/p)=I_p+\rho(p)/p.
\end{align*}
For $p>2$, the sieve factor gives
\begin{align*}
(1-\rho(p)/p)(I_p+\rho(p)/p)\leq I_p,
\end{align*}
since $I_p\geq1-\rho(p)/p$.
This comparison also applies at primes dividing the leading
coefficients or $h$.

At $p=2$, the theorem has no sieve factor.
Whether the local integrand vanishes depends only on parity.
If the sieve conditions exclude both classes, then
$I_2=0$ and every integer summand is zero.
Otherwise, the integrand is at least one on an admissible
class. It follows that $I_2\geq1/2$ and
\begin{align*}
I_2+\rho(2)/2\leq I_2+1\leq3I_2.
\end{align*}
This proves \eqref{eq:affine-rough}.

For \eqref{eq:affine-product}, the forms $\ell_1,\ell_2$ are
primitive and have resultant $-h$. Their coefficients are
bounded by $\exp(O_h(\sqrt L))$.
We apply \cite[Corollary~4.2]{Walsh-2026} with
$R=(x-r)/(AB)=x^{1+o(1)}$.
Since $-AB<r\leq0$, the summation interval is contained
in $(R,2R]$. The two cutoffs $R$ and $x/(AB)$ differ by
less than one.
Any extra local factor is uniformly bounded.
This gives \eqref{eq:affine-product}.
\end{proof}

We next apply this upper bound to products of dyadic intervals.
The resulting estimate will also be used in the short range
where the conditions $m>x^\gamma$ and $\bar m>x^\gamma$ differ.

\begin{lemma}\label{lem:factor-box}
The parameters $H_i$ satisfy $H_i\geq x^\gamma/2$
and $\prod_iH_i\asymp x$.
For each $1\leq i\leq j$, the coefficient $b_i$ is either
$a$ or $\widetilde a$. Then
\begin{equation}\label{eq:factor-box}
\sum_{m_i\asymp H_i}
Q_k(m_1\cdots m_j+h)\prod_{i=1}^jb_i(m_i)
\ll_{k,j,h,\gamma}xL^{k-j-1}.
\end{equation}
\end{lemma}

\begin{proof}
Each coefficient satisfies
$0\leq b_i(m)\leq(\gamma LV(Z))^{-1}$
and vanishes unless $P^-(m)>Z$.
We fix $m_2,\ldots,m_j$ on this support and set
$N=m_2\cdots m_j$.
Then $P^-(N)>Z$ and $N\leq x^{O(1)}$.
Lemma~\ref{lem:affine}, with the corrected local factors
at $p\mid N$, gives
\begin{equation}\label{eq:rough-affine-divisor}
\sum_{m\asymp H_1}\mathbf 1_{P^-(m)>Z}\tau_k(Nm+h)
\ll H_1V(Z)L^{k-1}.
\end{equation}

For $p\leq Z$, we have $p\nmid N$.
The sieve excludes the class $m\equiv0\;(\bmod p)$.
If also $p\nmid h$, the local factor is
\begin{align*}
(1-p^{-1})^{1-k}-p^{-1}
=1+\frac{k-2}{p}+O_k(p^{-2}).
\end{align*}
For $p>Z$ with $p\nmid N$, the factor is
$(1-p^{-1})^{1-k}$.
For $p\mid N$, we have $p\nmid h$ and thus $p\nmid Nm+h$
for every $m$. The local factor is $1$, so these primes do
not increase the product.
The finitely many factors at $p\mid h$ are absorbed into the
implied constant. Mertens' estimate bounds the product by
$V(Z)L^{k-1}$ and proves \eqref{eq:rough-affine-divisor}.

For the term containing $D_k$, we use
$D_k(Nm+h)\ll B_k\tau_k((Nm+h)_{\mathcal S})$.
The divisor weight now contributes $(\log y)^{k-1}$ in place
of $L^{k-1}$ in the same argument.
Since $B_k(\log y)^{k-1}\ll L^{k-1}$,
\eqref{eq:rough-affine-divisor} also holds with $\tau_k$
replaced by $D_k$.

In each remaining interval $m_i\asymp H_i$, an elementary
upper-bound sieve gives $O(H_iV(Z))$ integers with
$P^-(m_i)>Z$.
We combine these counts with the coefficient bounds
and \eqref{eq:rough-affine-divisor} to obtain
\eqref{eq:factor-box}.
\end{proof}
We first consider the range
$(m_1\cdots m_j)_{\mathcal S}\leq E$.
Then $0\leq\log(m_i/\bar m_i)\leq\sqrt L$ for each $i$.
If both $m_i$ and $\bar m_i$ exceed $x^\gamma$, we may replace
$1/\log m_i$ by $1/\log\bar m_i$ with relative error
$O_\gamma(L^{-1/2})$. If the indicators in
\eqref{eq:rough-factor} and \eqref{eq:prime-replacement-weight}
differ, then
\begin{equation}\label{eq:threshold}
 x^\gamma<m_i\leq x^\gamma E
\end{equation}
for at least one $i$.
On the common support, we expand the difference of the products
one factor at a time and obtain
\begin{align*}
 \Big|\prod_{i=1}^ja(m_i)-\prod_{i=1}^j\widetilde a(m_i)\Big|
 \ll_{j,\gamma}L^{-1/2}
 \sum_{\ell=1}^j\prod_{i<\ell}a(m_i)
                    \prod_{i\geq\ell}\widetilde a(m_i).
\end{align*}
There are $O(L^{j-1})$ products of dyadic intervals with
$\prod_iH_i\asymp x$. Lemma~\ref{lem:factor-box} and
$|\xi_k|\leq Q_k$ therefore bound the total denominator error by
\begin{align*}
 O\Big(L^{-1/2}L^{j-1}xL^{k-j-1}\Big)=O(xL^{k-5/2}).
\end{align*}
For any fixed coordinate $i$, at most
$O(L^{j-2}\log E)=O(L^{j-3/2})$ such products meet
\eqref{eq:threshold}. We apply Lemma~\ref{lem:factor-box} once more to bound the
error from the indicators by
\begin{align*}
 O\Big(L^{j-3/2}xL^{k-j-1}\Big)=O(xL^{k-5/2}).
\end{align*}
The bounded weight $\mathcal W$ changes only the constants in
these estimates. For $(m_1\cdots m_j)_{\mathcal S}>E$, we bound
the factorization multiplicities by $\tau_j(m_1\cdots m_j)$ and apply
Lemma~\ref{lem:moments}. Consequently, every $a(m_i)$ in
\eqref{eq:replacement} may be replaced by $\widetilde a(m_i)$,
with total error $o(xL^{k-2})$.

\subsection{Application of the graph estimate}

Proposition~\ref{prop:graph} requires sufficiently many prime sets
in which each endpoint has exactly one prime divisor.
 {We estimate these divisibility conditions by means of local
densities that incorporate the two divisor weights.}

We fix an integer $j\geq1$ and a set $ {\mathcal U}\subseteq\mathcal S$
whose primes all exceed $Z$. For sufficiently large $x$,
every prime in $ {\mathcal U}$ is coprime to $h$.
For each $p\in {\mathcal U}$, the normalizing factor is
\begin{align*}
I_{k,j,p}
=(1-p^{-1})^{1-k}+(1-p^{-1})^{1-j}-1.
\end{align*}
 {For nonnegative integers $a,b$, the corresponding local-density
weights are}
\begin{equation}\label{eq:record-density}
\mu_p(a,b)=\frac{1}{I_{k,j,p}}
\begin{cases}
1-2/p,
   & a=b=0,\\[3pt]
(1-p^{-1})\tau_k(p^a)p^{-a},
   & a\geq1,\ b=0,\\[3pt]
(1-p^{-1})\tau_j(p^b)p^{-b},
   & a=0,\ b\geq1,\\[3pt]
0,
   & a,b\geq1.
\end{cases}
\end{equation}
 {These local weights sum to one. We define their product on pairs} $(A,B)$ of $ {\mathcal U}$-smooth positive
integers by
\begin{align*}
\mu_{k,j}(A,B)
=\prod_{p\in {\mathcal U}}\mu_p(v_p(A),v_p(B)).
\end{align*}
For a bounded function of these pairs, we write
\begin{align*}
\langle F\rangle_{k,j}=
\sum_{\substack{A=A_{\mathcal U}\\B=B_{\mathcal U}}}
\mu_{k,j}(A,B)F(A,B).
\end{align*}
We adapt the local-density argument of
\cite[Lemma~4.3]{Walsh-2026} to include the sifting condition,
which contributes the factor $V(Z)^j$.

\begin{lemma}\label{lem:records}
For every nonnegative function $F(A,B)$ on pairs of
$ {\mathcal U}$-smooth positive integers with
$\|F\|_\infty\ll L^C$, where $C$ is fixed, we have
\begin{equation}\label{eq:rough-envelope}
\begin{aligned}
&\sum_{n\asymp x}
 Q_k(n)\tau_j(n-h)\mathbf 1_{P^-(n-h)>Z}
 F(n_{ {\mathcal U}},(n-h)_{ {\mathcal U}})\\
&\ll_{k,j,h}
 xL^{k+j-2}V(Z)^j {\langle F\rangle_{k,j}}
 +xL^{O(1)}\exp(-cL^{1/10}).
\end{aligned}
\end{equation}
The implied constant in the first term does not depend on
$F$ or $ {\mathcal U}$. The exponent of $L$ in the error
term may depend on $C$.
\end{lemma}

\begin{proof}
Lemma~\ref{lem:moments} allows us to restrict to
$n_{\mathcal S}\leq E$ and $(n-h)_{\mathcal S}\leq E$,
with the error in \eqref{eq:rough-envelope}.
We first consider the term containing $\tau_k(n)$.
Within this range, we fix $A=n_{ {\mathcal U}}$ and
$B=(n-h)_{ {\mathcal U}}$. These satisfy
$(A,B)=(AB,h)=1$ and $A,B\leq E$.

We use the representation $n=ABm+r$ from \eqref{eq:record-forms}.
The equalities $n_{ {\mathcal U}}=A$ and $(n-h)_{ {\mathcal U}}=B$
hold precisely when
$(\ell_1(m), {\mathcal U})=(\ell_2(m), {\mathcal U})=1$.
We apply \eqref{eq:affine-product} to the remaining divisor
weights, with the additional condition $P^-(n-h)>Z$.
At $p\in {\mathcal U}$, the local factor for the two residual forms
is $1-2/p$ if $p\nmid AB$, and $1-1/p$ if $p\mid AB$.
Together with $\tau_k(A)\tau_j(B)/(AB)$ and normalization
by $I_{k,j,p}$, these factors give {the local weights in}
\eqref{eq:record-density}.

For $p\leq Z$ with $p\nmid h$, the restriction $p\nmid n-h$
replaces $I_{k,j,p}$ by $J_p$, and
\begin{align*}
 \frac{J_p}{I_{k,j,p}}=1-\frac jp+O_{k,j}(p^{-2}).
\end{align*}
The factors at the finitely many primes dividing $h$ are
absorbed into the implied constant. The product over primes
outside $ {\mathcal U}$ is
\begin{align*}
 O_{k,j,h}\Big(
 \frac{L^{k+j-2}V(Z)^j}{\prod_{p\in {\mathcal U}}I_{k,j,p}}\Big).
\end{align*}
Together with the factors at $p\in {\mathcal U}$, this gives
\begin{align*}
\sum_{\substack{n\asymp x\\n_{ {\mathcal U}}=A,\ (n-h)_{ {\mathcal U}}=B}}
 \tau_k(n)\tau_j(n-h)\mathbf 1_{P^-(n-h)>Z}
 \ll xL^{k+j-2}V(Z)^j\mu_{k,j}(A,B).
\end{align*}
We then multiply by $F(A,B)$ and sum over $A,B\leq E$. Since $F$ is
nonnegative, we may extend the model sum to all pairs $(A,B)$,
which gives the required bound for $\tau_k(n)$.

For the term containing $D_k(n)$, we use
$D_k(n)\ll B_k\tau_k(n_{\mathcal S})$.
Before multiplication by $B_k$, the corresponding Euler product
outside $ {\mathcal U}$ is
\begin{align*}
 O_{k,j,h}\Big(
 \frac{(\log y)^{k-1}L^{j-1}V(Z)^j}
      {\prod_{p\in {\mathcal U}}I_{k,j,p}}\Big).
\end{align*}
Since $B_k(\log y)^{k-1}\asymp_k L^{k-1}$, this gives the same
upper bound. At every $p\in {\mathcal U}$,  {the local weights in}
\eqref{eq:record-density} are unchanged. We add the two estimates to obtain \eqref{eq:rough-envelope}.
\end{proof}

We now take $\mathcal U=\mathcal C$ in Lemma~\ref{lem:records}
and in the weighted mean $\langle\cdot\rangle_{k,j}$.
For $v=n-h$, a prime set satisfies the required condition when
\begin{equation}\label{eq:band-success}
\omega_g(n)=\omega_g(v)=1.
\end{equation}
We set $I_g(A,B)=\mathbf 1_{\omega_g(A)=\omega_g(B)=1}$
and $\beta_g=\langle I_g\rangle_{k,j}$.
For $p\in\mathcal P_g$, the local masses
\begin{align*}
U_p=\sum_{a\geq1}\mu_p(a,0)=\frac{k}{p}+O_{k,j}(p^{-2}),
\quad V_p=\sum_{b\geq1}\mu_p(0,b)=\frac{j}{p}+O_{k,j}(p^{-2})
\end{align*}
follow from \eqref{eq:record-density}. The product structure of the
weights gives
\begin{align*}
\beta_g
&=\prod_{p\in\mathcal P_g}(1-U_p-V_p)
 \sum_{\substack{p,q\in\mathcal P_g\\p\ne q}}
 \frac{U_pV_q}{(1-U_p-V_p)(1-U_q-V_q)}\\
&=\{kjV_g^2+o(1)\}\exp\{-(k+j)V_g+o(1)\}.
\end{align*}
Since $1/2\leq V_g\leq2$, we have
$\beta_g\geq\beta_{k,j}>0$ uniformly in $g$ for large $x$.

For one family of $N\asymp\log L$ disjoint sets, the
local weights factor, so an exponential weight gives
\begin{align*}
\Big\langle\mathbf 1_{\sum_g I_g<\beta_{k,j}N/2}\Big\rangle_{k,j}
&\leq e^{\beta_{k,j}N/2}
 \Big\langle e^{-\sum_g I_g}\Big\rangle_{k,j}=e^{\beta_{k,j}N/2}
 \prod_g\Big(1-(1-e^{-1})\beta_g\Big)\\
&\leq\exp\Big\{-\beta_{k,j}(1/2-e^{-1})N\Big\}
 \ll L^{-\eta}
\end{align*}
for some fixed $\eta>0$. Thus the pairs with too few sets satisfying \eqref{eq:band-success}
have weighted density bounded by a negative power of $L$. The $j$-fold Dirichlet convolution of $\widetilde a$ is
\begin{align*}
 \widetilde a^{*j}(v)
 =\sum_{m_1\cdots m_j=v}\prod_{i=1}^j\widetilde a(m_i).
\end{align*}
The coefficient bounds imply
\begin{equation}\label{eq:convolution-bound}
 \widetilde a^{*j}(v)
 \leq(\gamma LV(Z))^{-j}\tau_j(v)\mathbf 1_{P^-(v)>Z}.
\end{equation}
 {We apply the preceding estimate with an exponential weight to
each family in \eqref{eq:bands}, then add the two bounds for
the exceptional weighted densities.}
Equations \eqref{eq:rough-envelope} and \eqref{eq:convolution-bound}
now bound the total absolute contribution of {the tuples for
which either family has too few sets satisfying
\eqref{eq:band-success}} by $O_K(xL^{k-2-\eta})$.
We remove these tuples before expanding any signed weights.

For a family of $N$ sets,  {we write $I_g$ for
$I_g(n_{\mathcal C},v_{\mathcal C})$} and put
$r_0=\lceil\beta_{k,j}N/2\rceil$. The exact Boolean identity
\begin{equation}\label{eq:boolean}
 \mathbf 1_{\sum_gI_g\geq r_0}
 =\sum_{r=r_0}^N(-1)^{r-r_0}\binom{r-1}{r_0-1}
 \sum_{\substack{T\subseteq\{1,\ldots,N\}\\|T|=r}}
 \prod_{g\in T}I_g
\end{equation}
has total absolute coefficient mass at most $3^N$.
We apply this identity separately to the two families in
\eqref{eq:bands}.
Every subfamily in the resulting expansion contains
$\gg\log L$ sets from each family, and the total absolute
coefficient mass is $L^{O_{k,j}(1)}$.

We next remove the terms with $p^2\mid nv$ for some
$p\in\mathcal C$. Their weighted contribution is exponentially
small, even after multiplication by this coefficient mass.
We fix a subfamily indexed by $\mathcal G$, set $r=|\mathcal G|$
and $\mathcal P=\bigcup_{g\in\mathcal G}\mathcal P_g$, and impose
\eqref{eq:band-success} for every $g\in\mathcal G$.
Then \eqref{eq:invariance} gives $\xi_k(n)=k^r\alpha_{\mathcal P}(n)$.
 {There are $j^r$ ways to assign the $r$ distinct primes in
$\mathcal P$ dividing $v$ to the $j$ factors in the convolution.} The invariance of $\widetilde a$ consequently gives
\begin{equation}\label{eq:invariant-convolution}
\begin{aligned}
 \widetilde a^{*j}(v)&=j^rB(v),
 \quad B(v)=\widetilde a^{*j}(v_{\mathcal P'}),\quad |B(v)|&\leq(\gamma LV(Z))^{-j}\tau_j(v_{\mathcal P'}).
\end{aligned}
\end{equation}
Moreover, $B(pv)=B(v)$ for $p\in\mathcal P$.

The sets defining $\mathcal W$ are disjoint from $\mathcal C$,
so $\mathcal W(n)$ depends only on
$n_{\mathcal S\setminus\mathcal P}$. The definition of
$\Omega_{\mathcal P}$ also gives
\begin{align*}
 \prod_{g\in\mathcal G}I_g
 =\Big(\prod_{g\in\mathcal G}V_g^2\Big)
   \Omega_{\mathcal P}(n)\Omega_{\mathcal P}(v).
\end{align*}
Thus the expanded sum has the form required by
Proposition~\ref{prop:graph}, with $H=\mathcal W$ and $B$ as in
\eqref{eq:invariant-convolution}, apart from the factor
$(kj)^r\prod_{g\in\mathcal G}V_g^2$. This factor is bounded by
a fixed power of $L$.

Before applying that proposition, we restore the terms with
$p^2\mid nv$ for some $p\in\mathcal C$, since its sum does not
exclude them. The smallest prime in $\mathcal C$ exceeds
$\exp(L^{0.245})$. We count multiples of $p^2$ and apply the fixed divisor moment
bounds to obtain, for every fixed $C>0$,
\begin{align*}
 L^C\sum_{\substack{n\asymp x\\
 p^2\mid n(n-h)\ \text{for some }p\in\mathcal C}}
 Q_k(n)\tau_j(n-h)
 \ll xL^{C'}\exp(-cL^{0.245}).
\end{align*}
The same positive divisor bounds control the restored integrand,
including its fixed logarithmic weights. We may therefore restore these terms with an error of $o(xL^{-A})$
for every fixed $A>0$, even after we sum the absolute coefficients
in \eqref{eq:boolean}.
We choose the accuracy in Proposition~\ref{prop:graph} after
accounting for $(kj)^r\prod_{g\in\mathcal G}V_g^2$ and the
coefficient mass in \eqref{eq:boolean}. The resulting total is
$o(xL^{k-2})$. Together with \eqref{eq:replacement} and the
preceding error estimates, this proves \eqref{eq:composites}.

Finally, a squarefree integer with $j$ prime factors has exactly
$j!$ ordered factorizations into $j$ primes. Together with the
bound for repeated primes at the beginning of this section,
\eqref{eq:composites} gives
\begin{equation}\label{eq:composite-sum}
\begin{aligned}
 &\sum_{\substack{v\ \mathrm{composite}\\P^-(v)>x^\gamma}}
 \Phi((v+h)/x)\xi_k(v+h)\mathcal W(v+h)\\
 &\quad=\sum_{j=2}^{\lfloor2/\gamma\rfloor}\frac1{j!}
   \sum_{\substack{p_1,\ldots,p_j>x^\gamma\\p_i\ne p_\ell\ (i\ne\ell)}}
   \Phi((p_1\cdots p_j+h)/x)\xi_k(p_1\cdots p_j+h)
   \mathcal W(p_1\cdots p_j+h)+o(xL^{-A})\\
 &\quad=o(xL^{k-2}).
\end{aligned}
\end{equation}

\section{Proof of the main theorem}

We now combine the sieve estimate with the bound for composite
integers. This gives the required estimate over primes with an
auxiliary weight. We shall remove this weight by an average over
disjoint collections of prime sets, and then evaluate the mean of
the density model.

\subsection{An average over auxiliary weights}
We first need a positive upper bound that retains the mean of the
auxiliary weight. This bound will control the error when we replace
the normalized average of the weights by $1$.

\begin{lemma}\label{lem:prime-envelope}
The set $\mathcal U\subseteq\mathcal S$ consists of primes
exceeding $\exp(L^{1/10})$, none of which divides $h$.
The function $F$ is nonnegative and depends only on the part
supported on $\mathcal U$, with
$\|F\|_\infty\ll L^C$, where
$C$ is fixed. With $\mu_*$ as in Subsection~\ref{sec:marks}, we have
\begin{equation}\label{eq:prime-envelope}
\begin{aligned}
\sum_{p\asymp x}Q_k(p+h)F((p+h)_{ {\mathcal U}})\ll_{k,h}xL^{k-2} {\sum_{a=a_{\mathcal S}}\mu_*(a)F(a_{\mathcal U})}
  +xL^{O(1)}\exp(-cL^{1/10}).
\end{aligned}
\end{equation}
The constant in the first term is independent of $F$ and $ {\mathcal U}$.
The exponent of $L$ in the second term may depend on $C$.
\end{lemma}

\begin{proof}
With $n=p+h$, Lemma~\ref{lem:moments} allows us to restrict to
$n_{\mathcal S}\le E$, with the error stated in
\eqref{eq:prime-envelope}. For sufficiently large $x$, every prime in
$ {\mathcal U}$ is less than $x^{1/10}$. We bound the prime indicator by
$\mathbf 1_{P^-(n-h)>x^{1/10}}$.

For each fixed small-prime part $A=n_{ {\mathcal U}}\le E$, we apply
\eqref{eq:affine-product} with $B=1$. At a prime $p\in {\mathcal U}$,
the resulting local masses are
\begin{equation}\label{eq:prime-density}
\begin{cases}
1-2/p,&v_p(A)=0,\\
(1-p^{-1})\tau_k(p^a)p^{-a},&v_p(A)=a\ge1.
\end{cases}
\end{equation}
Their sum is $J_p$. We divide each mass by $J_p$ to recover the
 {corresponding normalized local-density weight.} The product of the remaining local factors is
\begin{align*}
O_{k,h}\Big(\frac{L^{k-2}}{\prod_{p\in {\mathcal U}}J_p}\Big).
\end{align*}
Consequently, the part containing $\tau_k(n)$ satisfies
\begin{align*}
\sum_{\substack{n\asymp x\\n_{ {\mathcal U}}=A}}
\tau_k(n)\mathbf 1_{P^-(n-h)>x^{1/10}}
\ll_{k,h}xL^{k-2} {\sum_{\substack{a=a_{\mathcal S}\\a_{\mathcal U}=A}}\mu_*(a)}.
\end{align*}
We multiply by $F(A)$ and sum over $A\leq E$.
Positivity then gives the weighted mean
$\sum_{a=a_{\mathcal S}}\mu_*(a)F(a_{\mathcal U})$.

For the contribution of $D_k(n)$, we use the bound
$D_k(n)\ll_k B_k\tau_k(n_{\mathcal S})$.
Before multiplication by $B_k$, the Euler product over primes
outside $ {\mathcal U}$ is
\begin{align*}
O_{k,h}\Big(\frac{(\log y)^{k-1}}
 {L\prod_{p\in {\mathcal U}}J_p}\Big).
\end{align*}
As in the proof of Lemma~\ref{lem:records}, the relation
$B_k(\log y)^{k-1}\asymp_k L^{k-1}$ gives the same bound.
The local masses in \eqref{eq:prime-density} remain unchanged.
The contribution from $n_{\mathcal S}>E$ is included in the
error, which proves \eqref{eq:prime-envelope}.
\end{proof}

We now form an average of weights from disjoint collections
of auxiliary prime sets. We compare this average with $1$
on the primes.

\begin{proposition}\label{prop:prime-centering}
For fixed $k\ge2$, $h\ne0$ and nonnegative
$\Phi\in C_c^\infty((0,\infty))$, we have
\begin{equation}\label{eq:prime-centering}
\sum_p\Phi((p+h)/x)\xi_k(p+h)=o(xL^{k-2}).
\end{equation}
\end{proposition}

\begin{proof}
We choose a sufficiently large even integer $s$ and set
$\gamma=c_0/(8s+8)$. We then choose $K$ as in
Section~\ref{sec:composites}. Every integer $v$ on the support
with $P^-(v)>x^\gamma$ is either prime or a composite product
covered by \eqref{eq:composites}. Proposition~\ref{prop:presieve}
and \eqref{eq:composite-sum} therefore give, for each fixed
collection of sets defining $\mathcal W$,
\begin{equation}\label{eq:marked-primes}
\limsup_{x\to\infty}
\frac{\Big|\sum_p\Phi((p+h)/x)\xi_k(p+h)\mathcal W(p+h)\Big|}
     {xL^{k-2}m_{\mathcal W}}
\ll_{k,h,\Phi}\frac{e^{-s}}{\gamma}.
\end{equation}
The implied constant depends neither on $K$ nor on the
exponents defining the prime sets.

For fixed $b\ge1$, we choose $bK$ pairwise disjoint prime sets
and divide them into $b$ collections of $K$ sets.
The corresponding weights are $\mathcal W_1,\ldots,\mathcal W_b$,
with $m_i= {\langle\mathcal W_i\rangle_*}$.
Their normalized average is
\begin{equation}\label{eq:mark-average}
Z_b(n)=\frac1b\sum_{i=1}^{b}\frac{\mathcal W_i(n)}{m_i}.
\end{equation}
For $i\neq j$, the weights $\mathcal W_i$ and $\mathcal W_j$
depend on disjoint prime sets, {so the local-density weights factor. Hence}
\begin{align*}
 {\langle\mathcal W_i\mathcal W_j\rangle_*}=m_i m_j.
\end{align*}
By \eqref{eq:mark-moments}, the diagonal terms satisfy
$ {\langle(\mathcal W_i/m_i)^2\rangle_*}\leq C_{K,k}$.
Since $ {\langle Z_b\rangle_*}=1$, we may expand the square to obtain
\begin{equation}\label{eq:mark-square}
 {\langle|Z_b-1|^2\rangle_*}
=\frac1{b^2}\sum_{i=1}^b
\Big\{ {\langle(\mathcal W_i/m_i)^2\rangle_*}-1\Big\}
\leq\frac{C_{K,k}}b.
\end{equation}
For fixed $b$ and $K$, every prime in the union of these sets exceeds
$\exp(L^{1/10})$, so Lemma~\ref{lem:prime-envelope} applies to
the part supported on their union.

We apply that lemma with $F=1$ and $F=|Z_b-1|^2$. Since
$|\xi_k|\le Q_k$, the weighted Cauchy--Schwarz inequality and
\eqref{eq:mark-square} give
\begin{equation}\label{eq:unmarking}
\limsup_{x\to\infty}
\frac{\Big|\sum_p\Phi((p+h)/x)\xi_k(p+h)(Z_b(p+h)-1)\Big|}
     {xL^{k-2}}
\ll_{k,h,\Phi}\sqrt{\frac{C_{K,k}}b}.
\end{equation}
The bound for $|Z_b-1|^2$ may depend on $b$ and $K$. In
\eqref{eq:prime-envelope}, this dependence affects only an error
that tends to zero as $x\to\infty$ with $b$ and $K$ fixed.

We average \eqref{eq:marked-primes} over the weights in
\eqref{eq:mark-average} and use \eqref{eq:unmarking} to obtain
\begin{align*}
\limsup_{x\to\infty}
\frac{\Big|\sum_p\Phi((p+h)/x)\xi_k(p+h)\Big|}{xL^{k-2}}
\le C_{k,h,\Phi}\frac{e^{-s}}{\gamma}
 +C'_{k,h,\Phi}\sqrt{\frac{C_{K,k}}b}.
\end{align*}
For each fixed $s$, we first take the limit superior as $x\to\infty$
and then let $b\to\infty$. Finally, we let $s\to\infty$. Since
$\gamma=c_0/(8s+8)$, we have $e^{-s}/\gamma\to0$.
This proves \eqref{eq:prime-centering}.
\end{proof}

\begin{remark}
We choose the parameters in the following order: $s$, $\gamma$,
the moment truncations and Type II accuracy, $K$, and finally $b$
and the exponents defining the prime sets. We then let $x\to\infty$,
followed by $b\to\infty$ and $s\to\infty$. No estimate uniform
in $K$ or $b$ growing with $x$ is required.
\end{remark}

\subsection{The model mean and completion of the proof}

Proposition~\ref{prop:prime-centering} reduces the prime correlation
to the mean of the density model. The model involves only small smooth divisors. We can therefore
evaluate its mean by the Bombieri--Vinogradov theorem.

\begin{proposition}\label{prop:model-mean}
For fixed $k\ge2$, $h\ne0$ and
$\Phi\in C_c^\infty((0,\infty))$, we have
\begin{equation}\label{eq:model-prime-mean}
\sum_n\Lambda(n)\Phi((n+h)/x)D_k(n+h)
=C_{k,h}xL^{k-1}\int_0^\infty\Phi(t)\,dt+o(xL^{k-1}).
\end{equation}
\end{proposition}

\begin{proof}
We first evaluate the sum in \eqref{eq:model-prime-mean} with
$D_k$ replaced by $D_k^0$. The identity
\begin{align*}
\tau_k((n+h)_{\mathcal S})
=\sum_{\substack{d\mid n+h\\d=d_{\mathcal S}}}\tau_{k-1}(d)
\end{align*}
expresses the model in terms of smooth divisors. If a divisor
$d>E$ occurs on the right, then $(n+h)_{\mathcal S}>E$.
Lemma~\ref{lem:moments}, together with $\Lambda(n)\le\log n$,
therefore allows us to restrict to $d\le E$ with a negligible error.

If $(d,h)>1$, $d\mid n+h$ and $\Lambda(n)\ne0$, then $n$ is
a power of a prime dividing $h$. There are $O_h(L)$ such prime
powers on the support. Their contribution is negligible by the
pointwise bound $\tau_k(n+h)\ll_{k,\varepsilon}x^\varepsilon$.
For $(d,h)=1$, we define
\begin{align*}
R_d=\sum_{n\equiv-h\bmod d}\Lambda(n)\Phi((n+h)/x)
 -\frac{x}{\varphi(d)}\int_0^\infty\Phi(t)\,dt.
\end{align*}
By partial summation and the Bombieri--Vinogradov
theorem, we obtain
\begin{equation}\label{eq:bv}
\sum_{\substack{d\le E\\(d,h)=1}}|R_d|\ll_A xL^{-A}
\end{equation}
for every fixed $A>0$. We also have $|R_d|\ll xL/d$. We apply the Cauchy--Schwarz inequality to obtain
\begin{align*}
\sum_{\substack{d\le E\\(d,h)=1}}\tau_{k-1}(d)|R_d|\ll\Big(\sum_{d\le E}\frac{\tau_{k-1}(d)^2}{d}\Big)^{1/2}
 \Big(xL\sum_{\substack{d\le E\\(d,h)=1}}|R_d|\Big)^{1/2}\ll_A xL^{-A}.
\end{align*}
The first factor grows by at most a fixed power of $L$, which we
absorb by choosing a larger accuracy in \eqref{eq:bv}.

We next extend the sum of the main terms to all $d=d_{\mathcal S}$
with $(d,h)=1$. Rankin's trick with exponent $L^{-2/5}$
bounds the tail by $L^{O_k(1)}\exp(-L^{1/10})$. Thus
\begin{equation}\label{eq:model-mean}
\begin{split}
\sum_n\Lambda(n)\Phi((n+h)/x)D_k^0(n+h)
={}&xB_k\Big(\int_0^\infty\Phi(t)\,dt\Big)
\prod_{\substack{p\le y\\p\nmid h}}
\Big\{1+\frac{(1-p^{-1})^{-(k-1)}-1}{1-p^{-1}}\Big\}\\
&+o(xL^{k-1}).
\end{split}
\end{equation}
For $p\nmid h$, we combine this local factor with
$(1-p^{-1})^{k-1}$ from $B_k$ and obtain
\begin{align*}
(1-p^{-1})^{k-1}
\Big\{1+\frac{(1-p^{-1})^{-(k-1)}-1}{1-p^{-1}}\Big\}
=\frac{1}{1-p^{-1}}-\frac{(1-p^{-1})^{k-2}}p.
\end{align*}
For $p\mid h$, the factor remains $(1-p^{-1})^{k-1}$.
By the absolute convergence of the product in \eqref{eq:constant},
the main term in \eqref{eq:model-mean} is therefore the one in
\eqref{eq:model-prime-mean}.

Finally, we apply Lemma~\ref{lem:model} on
$(n+h)_{\mathcal S}\le E$ to replace $D_k^0(n+h)$ by
$D_k(n+h)$ with relative error $O_k(L^{-1/2})$. We choose a
nonnegative $\Phi_0\in C_c^\infty((0,\infty))$ with
$\Phi_0\ge|\Phi|$ and apply \eqref{eq:model-mean} with
$\Phi_0$ in place of $\Phi$.
The resulting error is $O(xL^{k-3/2})$. The part with
$(n+h)_{\mathcal S}>E$ is negligible by Lemma~\ref{lem:moments}.
This proves \eqref{eq:model-prime-mean}.
\end{proof}

\begin{proof}[Proof of Theorem~\ref{thm:main}]
We first consider a fixed nonnegative smooth weight $\Phi$. By
\eqref{eq:prime-envelope} with $F=1$,
\begin{align*}
\sum_{p\asymp x}|\xi_k(p+h)|\ll xL^{k-2}.
\end{align*}
On the support of $\Phi((p+h)/x)$, we have
$\log p=L+O_\Phi(1)$. Proposition~\ref{prop:prime-centering}
therefore implies
\begin{equation}\label{eq:prime-log}
\sum_p(\log p)\Phi((p+h)/x)\xi_k(p+h)=o(xL^{k-1}).
\end{equation}

Prime powers of exponent at least two contribute
$O_{k,h,\varepsilon}(x^{1/2+\varepsilon}L^{O_k(1)})$, by the
pointwise divisor bound and the corresponding bound for $D_k$.
Thus \eqref{eq:prime-log} also holds with the sum over all positive
integers and weight $\Lambda(n)$. We combine this estimate with Proposition~\ref{prop:model-mean}
and the identity $\tau_k=D_k+\xi_k$ to obtain
\begin{align*}
\sum_n\Lambda(n)\Phi((n+h)/x)\tau_k(n+h)
=C_{k,h}xL^{k-1}\int_0^\infty\Phi(t)\,dt+o(xL^{k-1}).
\end{align*}

To remove the smooth weight, we fix $0<\varepsilon<1/4$.
We choose nonnegative smooth functions $\Phi_-$ and $\Phi_+$ such that
$\Phi_-\leq1$, $\Phi_-$ is supported in $(1+\varepsilon,2-\varepsilon)$,
and $\Phi_+\geq1$ on $[1-\varepsilon,2+\varepsilon]$, with
\begin{align*}
\int_0^\infty\Phi_-(t)\,dt=1+O(\varepsilon),\quad
\int_0^\infty\Phi_+(t)\,dt=1+O(\varepsilon).
\end{align*}
For sufficiently large $x$, we have $|h|/x<\varepsilon/2$.
The nonnegativity of $\Lambda(n)\tau_k(n+h)$ then gives
\begin{align*}
\sum_n\Lambda(n)\tau_k(n+h)\Phi_-((n+h)/x)
&\leq\sum_{x<n\leq2x}\Lambda(n)\tau_k(n+h)\\
&\leq\sum_n\Lambda(n)\tau_k(n+h)\Phi_+((n+h)/x).
\end{align*}
We apply the smooth asymptotic to both bounds, first let
$x\to\infty$, and then let $\varepsilon\to0$.
This proves Theorem \ref{thm:main}.
\end{proof}

\section*{Acknowledgements}
The authors would like to thank ChatGPT for its assistance in understanding the methods in \cite{Walsh-2026}. The authors also thank ChatGPT for helping to check the manuscript.

\bibliographystyle{amsplain}
\bibliography{titchmarsh_references}

@article{Assing-Blomer-Li-2021,
  author        = {Assing, E. and Blomer, V. and Li, J.},
  title         = {{Uniform Titchmarsh divisor problems}},
  journal       = {Adv. Math.},
  volume        = {393},
  year          = {2021},
  pages         = {108076},
  doi           = {10.1016/j.aim.2021.108076},
  eprint        = {2005.13915v1},
  archivePrefix = {arXiv},
  primaryClass  = {math.NT},
  url           = {https://arxiv.org/abs/2005.13915v1}
}

@article{Bombieri-Friedlander-Iwaniec-1986,
  author  = {Bombieri, E. and Friedlander, J. B. and Iwaniec, H.},
  title   = {{Primes in arithmetic progressions to large moduli}},
  journal = {Acta Math.},
  volume  = {156},
  year    = {1986},
  pages   = {203--251},
  doi     = {10.1007/BF02399204}
}

@article{Deshouillers-Iwaniec-1982,
  author  = {Deshouillers, J. and Iwaniec, H.},
  title   = {{Kloosterman sums and Fourier coefficients of cusp forms}},
  journal = {Invent. Math.},
  volume  = {70},
  number  = {2},
  year    = {1982},
  pages   = {219--288},
  doi     = {10.1007/BF01390728}
}

@article{Drappeau-2017,
  author  = {Drappeau, S.},
  title   = {{Sums of Kloosterman sums in arithmetic progressions, and the error term in the dispersion method}},
  journal = {Proc. Lond. Math. Soc. (3)},
  volume  = {114},
  number  = {4},
  year    = {2017},
  pages   = {684--732}
}

@article{Drappeau-Topacogullari-2019,
  author  = {Drappeau, S. and Topacogullari, B.},
  title   = {{Combinatorial identities and Titchmarsh's divisor problem for multiplicative functions}},
  journal = {Algebra Number Theory},
  volume  = {13},
  number  = {10},
  year    = {2019},
  pages   = {2383--2425}
}

@article{Fiorilli-2012,
  author  = {Fiorilli, D.},
  title   = {{Residue classes containing an unexpected number of primes}},
  journal = {Duke Math. J.},
  volume  = {161},
  number  = {15},
  year    = {2012},
  pages   = {2923--2943}
}

@article{Fouvry-1985,
  author  = {Fouvry, {\'E}.},
  title   = {{Sur le probl\`eme des diviseurs de Titchmarsh}},
  journal = {J. reine angew. Math.},
  volume  = {357},
  year    = {1985},
  pages   = {51--76},
  doi     = {10.1515/crll.1985.357.51}
}

@article{Halberstam-1967,
  author  = {Halberstam, H.},
  title   = {{Footnote to the Titchmarsh--Linnik divisor problem}},
  journal = {Proc. Amer. Math. Soc.},
  volume  = {18},
  year    = {1967},
  pages   = {187--188},
  doi     = {10.1090/S0002-9939-1967-0204379-6}
}

@article{Henriot-2012,
  author  = {Henriot, K.},
  title   = {{Nair--Tenenbaum bounds uniform with respect to the discriminant}},
  journal = {Math. Proc. Cambridge Philos. Soc.},
  volume  = {152},
  number  = {3},
  year    = {2012},
  pages   = {405--424}
}

@article{Henriot-2014,
  author  = {Henriot, K.},
  title   = {{Nair--Tenenbaum uniform with respect to the discriminant---Erratum}},
  journal = {Math. Proc. Cambridge Philos. Soc.},
  volume  = {157},
  number  = {2},
  year    = {2014},
  pages   = {375--377}
}

@book{Iwaniec-Kowalski-2004,
  author    = {Iwaniec, H. and Kowalski, E.},
  title     = {{Analytic number theory}},
  series    = {American Mathematical Society Colloquium Publications},
  volume    = {53},
  publisher = {American Mathematical Society},
  address   = {Providence, RI},
  year      = {2004}
}

@article{Linnik-1960,
  author  = {Linnik, Yu. V.},
  title   = {{An asymptotic formula in an additive problem of Hardy--Littlewood}},
  journal = {Izv. Akad. Nauk SSSR Ser. Mat.},
  volume  = {24},
  number  = {5},
  year    = {1960},
  pages   = {629--706},
  note    = {In Russian},
  url     = {https://www.mathnet.ru/eng/im3674}
}

@book{Linnik-1963,
  author    = {Linnik, Yu. V.},
  title     = {{The dispersion method in binary additive problems}},
  series    = {Translations of Mathematical Monographs},
  volume    = {4},
  publisher = {American Mathematical Society},
  address   = {Providence, RI},
  year      = {1963}
}

@misc{Lu-Wang-Wu-2026,
  author        = {L{\"u}, G. S. and Wang, Z. W. and Wu, J.},
  title         = {{On the Hardy--Littlewood--Chowla conjecture}},
  howpublished  = {arXiv:2610.20242v1},
  year          = {2026},
  eprint        = {2610.20242v1},
  archivePrefix = {arXiv},
  primaryClass  = {math.NT},
  url           = {https://arxiv.org/abs/2610.20242v1}
}

@article{Matomaki-Radziwill-Tao-2019a,
  author  = {Matom{\"a}ki, K. and Radziwi{\l}{\l}, M. and Tao, T.},
  title   = {{Correlations of the von Mangoldt and higher divisor functions I. Long shift ranges}},
  journal = {Proc. Lond. Math. Soc. (3)},
  volume  = {118},
  number  = {2},
  year    = {2019},
  pages   = {284--350}
}

@article{Matomaki-Radziwill-Tao-2019b,
  author  = {Matom{\"a}ki, K. and Radziwi{\l}{\l}, M. and Tao, T.},
  title   = {{Correlations of the von Mangoldt and higher divisor functions II. Divisor correlations in short ranges}},
  journal = {Math. Ann.},
  volume  = {374},
  number  = {1--2},
  year    = {2019},
  pages   = {793--840}
}

@article{Matomaki-Radziwill-Shao-Tao-Teravainen-2026,
  author  = {Matom{\"a}ki, K. and Radziwi{\l}{\l}, M. and Shao, X. C. and Tao, T. and Ter{\"a}v{\"a}inen, J.},
  title   = {{Higher uniformity of arithmetic functions in short intervals II. Almost all intervals}},
  journal = {Invent. Math.},
  volume  = {244},
  year    = {2026},
  pages   = {967--1091}
}

@misc{OpenAI-2026c,
  author       = {{OpenAI}},
  title        = {{The Poisson--Dirichlet law for prime predecessors}},
  howpublished = {OpenAI Math Release preprint},
  year         = {2026},
  month        = {24 September},
  url          = {https://github.com/openai/math/tree/main/preprints/The-Poisson-Dirichlet-Law-for-Prime-Predecessors-September-24-2026}
}

@misc{OpenAI-2026f,
  author       = {{OpenAI}},
  title        = {{Weighted dilation graphs, smooth shifted primes and totient fibers}},
  howpublished = {OpenAI Math Release preprint},
  year         = {2026},
  month        = {24 September},
  url          = {https://github.com/openai/math/tree/main/preprints/Weighted-Dilation-Graphs-Smooth-Shifted-Primes-and-Totient-Fibers-September-24-2026}
}

@article{Rodriquez-1965,
  author  = {Rodriquez, G.},
  title   = {{Sul problema dei divisori di Titchmarsh}},
  journal = {Boll. Un. Mat. Ital. (3)},
  volume  = {20},
  number  = {3},
  year    = {1965},
  pages   = {358--366},
  url     = {https://www.bdim.eu/item?id=BUMI_1965_3_20_3_358_0}
}

@article{Titchmarsh-1930,
  author  = {Titchmarsh, E. C.},
  title   = {{A divisor problem}},
  journal = {Rend. Circ. Mat. Palermo},
  volume  = {54},
  year    = {1930},
  pages   = {414--429},
  doi     = {10.1007/BF03021203},
  note    = {Correction, ibid. 57 (1933), 478--479}
}

@misc{Walsh-2026,
  author        = {Walsh, M. N.},
  title         = {{Shifted convolutions of divisor functions}},
  howpublished  = {arXiv:2610.40439v1},
  year          = {2026},
  eprint        = {2610.40439v1},
  archivePrefix = {arXiv},
  primaryClass  = {math.NT},
  url           = {https://arxiv.org/abs/2610.40439v1}
}
\end{document}